\documentclass[10pt]{article}
\usepackage{amsmath,amsfonts,geompsfi}

\renewcommand{\bold}[1]{\medskip \noindent {\bf \boldmath #1
                        }\nopagebreak[4]}

\newcommand{\qed}{\nopagebreak[4]\hspace{.2cm} $\square$ \pagebreak[2]\medskip}
\newtheorem{theorem}{Theorem}[section]

\newcommand{\cx}{{\mathbb C}}
\newcommand{\half}{{\mathbb H}}
\newcommand{\integers}{{\mathbb Z}}

\newcommand{\reals}{{\mathbb R}}

\newcommand{\tube}{{\mathbb T}}
\newcommand{\cusp}{{\mathbb P}}
\newcommand{\makefig}[3]{
	\begin{figure}[htbp]
        \refstepcounter{figure}
	\label{#2}
        \begin{center}~
		#3~\\
		\medskip
                {\sf Figure \thefigure.  #1}
        \end{center}
	\medskip
	\end{figure}
}

\newenvironment{pf*}[1]{%
 \begin{proof}[#1]%
}{ 
 \end{proof}
}

\newcommand{\ital}[1]{\medskip \noindent {\em #1 }\nopagebreak[4]}

\newcommand{\bdry}{\partial}

\newcommand{\closure}{\overline}
\newcommand{\compos}{\circ}

\newcommand{\del}{\partial}
\newcommand{\dirsum}{\oplus}

\newcommand{\disjunion}{\sqcup}

\usepackage{amssymb}
\newcommand{\nullset}{\varnothing}

\newcommand{\wt}{\widetilde}

\newcommand{\chat}{\widehat{\cx}}

\newcommand{\zed}{\integers}

\newcommand{\Aut}{\mbox{\rm Aut}}

\newcommand{\core}{\mbox{\rm core}}

\newcommand{\diam}{\mbox{\rm diam}}

\newcommand{\hull}{\mbox{\rm hull}}

\newcommand{\inj}{\mbox{\rm inj}}
\newcommand{\id}{\mbox{\rm id}}
\newcommand{\interior}{\mbox{\rm int}}

\newcommand{\Isom}{\mbox{\rm Isom}}
\newcommand{\join}{\mbox{\rm join}}

\newcommand{\PSL}{\mbox{\rm PSL}}

\newcommand{\Teich}{\mbox{\rm Teich}}

\newcommand{\vol}{\mbox{\rm vol}}

\newtheorem{prop}[theorem]{Proposition}
\newtheorem{lem}[theorem]{Lemma}
\newtheorem{cor}[theorem]{Corollary}
\newtheorem{conj}[theorem]{Conjecture}
\newtheorem{defn}[theorem]{Definition}

\newcommand{\calA}{{\mathcal A}}

\newcommand{\calC}{{\mathcal C}}

\newcommand{\calE}{{\mathcal E}}

\newcommand{\calL}{{\mathcal L}}
\newcommand{\calM}{{\mathcal M}}
\newcommand{\calN}{{\mathcal N}}

\newcommand{\calP}{{\mathcal P}}
\newcommand{\calQ}{{\mathcal Q}}

\newcommand{\calS}{{\mathcal S}}

\newcommand{\ml}{{\calM \calL}}
\newcommand{\pl}{{\calP \calL}}
\newcommand{\ps}{{\calP \calS}}
\newcommand{\el}{{\calE \calL}}

\usepackage{euscript}

\newcommand{\eS}{{\EuScript S}}

\date{December 12, 2002}
\title{\vspace{-.5in}
{\bf On the density of geometrically finite Kleinian groups}
\vspace{.2in}}
\author{Jeffrey F. Brock\thanks{Research supported by an NSF
Postdoctoral Fellowship and NSF research grants.}
\ and Kenneth W. Bromberg\thanks{Research supported by
NSF research grants and the Clay Mathematics Institute. \hfill
}
}

\newcommand{\short}{{\bf short}}

\begin{document}

\maketitle

\begin{abstract}
\noindent The density conjecture of Bers, Sullivan and Thurston predicts that
each complete hyperbolic 3-manifold $M$ with finitely generated
fundamental group is an algebraic limit of geometrically finite
hyperbolic 3-manifolds.  We prove that the conjecture obtains for each
complete hyperbolic 3-manifold with no cusps and incompressible ends.

\end{abstract}
{\scriptsize
\tableofcontents
}
\section{Introduction}

In the 1970s, work of W. P. Thurston revolutionized the study of
Kleinian groups and their 3-dimensional hyperbolic quotients.
Nevertheless, a complete topological and geometric classification of
hyperbolic 3-manifolds persists as a fundamental unsolved problem.

Even for {\em tame} hyperbolic 3-manifolds $N = \half^3 / \Gamma$,
where $N$ has tractable topology ($N$ is homeomorphic to the interior
of a compact 3-manifold), the correct picture of the range of
complete hyperbolic structures on $N$ remains conjectural.
\smallskip

On the other hand, {\em geometrically finite} hyperbolic 3-manifolds
are completely parameterized by an elegant deformation theory.  As an
approach to a general classification, Thurston proposed a program to
extend this parameterization to all hyperbolic 3-manifolds with
finitely generated fundamental group \cite{Thurston:survey}.  A
critical, and as yet 
unyielding obstacle is the {\em density conjecture}:
\begin{conj}[Bers-Sullivan-Thurston]
Let $M$ be a complete hyperbolic 3-manifold with finitely generated
fundamental group.  Then $M$ is a limit of geometrically finite hyperbolic
3-manifolds.
\label{conjecture:density}
\end{conj}
Our main result is the following.
\begin{theorem}
Let $M$ be a complete hyperbolic 3-manifold with finitely generated
fundamental group, incompressible ends, and no cusps.  Then $M$ is an
algebraic limit of geometrically finite hyperbolic 3-manifolds.
\label{theorem:main}
\end{theorem}
We call $M$ {\em geometrically finite} if its {\em convex core}, the
minimal convex subset of $M$, has finite volume.  The manifold $M$ is
the quotient of $\half^3$ by a {\em Kleinian group} $\Gamma$, a
discrete, torsion free subgroup of the orientation preserving
isometries of hyperbolic 3-space.  Then $M =
\half^3/\Gamma$ is an {\em algebraic limit} of $M_i =
\half^3/\Gamma_i$ if there are isomorphisms $\rho_i
\colon
\Gamma \to \Gamma_i$ so that after conjugating the groups $\Gamma_i$
in $\Isom^+\half^3$ if necessary, we have $\rho_i(\gamma) \to
\gamma$ for each $\gamma \in \Gamma$.  We say $M$ has {\em
incompressible ends} if it is homotopy equivalent to a compact
submanifold with incompressible boundary.

Marden and Sullivan proved that the interior of the deformation space
of complete hyperbolic 3-manifolds homotopy equivalent to $M$ consists
of such geometrically finite hyperbolic 3-manifolds
(see \cite{Marden:kgs,Sullivan:QCDII}).  Then
Conjecture~\ref{conjecture:density} predicts that the deformation
space is the closure of its interior.

\smallskip

Theorem~\ref{theorem:main}
generalizes the recent result of the second author
\cite{Bromberg:bers}, which applies to cusp-free {\em singly
degenerate} manifolds 
$M$ with the homotopy type of a surface.  In
that case, the result gives a partial solution to an earlier version
of Conjecture~\ref{conjecture:density} formulated by L. Bers in
\cite{Bers:bdry}. 
For the modern formulation, see 
\cite{Sullivan:QCDII}, and
\cite{Thurston:survey}.

Fundamental in the treatment of each case is the use of 3-dimensional
hyperbolic cone-manifolds, namely, 3-manifolds that are hyperbolic
away from a closed geodesic cone-type singularity.  The theory of
deformations of these manifolds that change only the cone-angle,
developed by C. Hodgson, S. Kerckhoff
\cite{Hodgson:Kerckhoff:rigidity,
Hodgson:Kerckhoff:bounds,
Hodgson:Kerckhoff:shape} and the second author
\cite{Bromberg:thesis,Bromberg:Schwarzian},
is instrumental in our study.
In particular, the recent innovations of
\cite{Hodgson:Kerckhoff:bounds} and \cite{Hodgson:Kerckhoff:shape}
have extended the theory to treat the setting of arbitrary
cone-angles, whereas \cite{Hodgson:Kerckhoff:rigidity} treats only the
case of cone-angle at most $2 \pi$ (see
\cite{Hodgson:Kerckhoff:harmonic} for an expository account).  These
estimates are essential to results of \cite{Bromberg:Schwarzian} and
their generalizations here.

Though the power of cone-deformations has been amply demonstrated in
the proof of the orbifold theorem and the study of hyperbolic
Dehn-surgery space developed by Hodgson and Kerckhoff, we hope the
present study will suggest its wider applicability as a new tool in
the study of deformation spaces of infinite volume hyperbolic
3-manifolds.

The principal application of the cone-deformation theory here is its
ability to control the geometric effect of a cone-deformation that
decreases the cone-angle at a short cone-singularity.  Since each
simple closed geodesic in a hyperbolic 3-manifold may be regarded as a
cone ``singularity'' with cone-angle $2 \pi$, we obtain control on how
a geometrically finite structure with a short closed geodesic differs
from the complete hyperbolic structure on the manifold with the same
conformal boundary and the short geodesic removed (the resulting cusp
may be viewed as a cone-singularity with cone-angle $0$).

A central result of the paper is a {\em drilling theorem},
giving an example of this type of control.  Here is a version
applicable to complete, smooth hyperbolic structures:
\begin{theorem}{\sc The Drilling Theorem}
Let $M$ be a geometrically finite hyperbolic 3-manifold.  For each $L
>1$, there is an $\ell >0$ so that if $c$ is a geodesic in $M$ with
length $\ell_M(c) < \ell$, there is an $L$-bi-Lipschitz
diffeomorphism of pairs
$$h \colon (M \setminus \tube(c), \bdry \tube(c))  \to (M_0 \setminus
\cusp(c), \bdry \cusp(c))$$
where $M\setminus \tube(c)$ denotes
 the complement of a standard tubular neighborhood of $c$ in $M$,
$M_0$ denotes the complete hyperbolic structure on $M \setminus c$,
and $\cusp(c)$ denotes a standard rank-2 cusp corresponding to $c$.
\label{theorem:smooth:drill}
\end{theorem}
(See Theorem~\ref{theorem:bilip} for a more precise version).

The drilling theorem and its algebraic antecedents in
\cite{Bromberg:Schwarzian} are 
reminiscent of the essential estimates needed to control the algebraic
effect of other types of pinching deformations.  Such estimates have
been used to show (for example) the density of maximal cusps in
boundaries of deformation spaces
\cite{McMullen:cusps,CCHS:density}.  While these estimates give
algebraic control over pinching short curves in the conformal
boundary, a very short geodesic in $M$ can have large length on the
conformal boundary of $M$. The drilling theorem, by contrast, applies
to any short geodesic in $M$.

The drilling theorem has proven to be of general use in the study of
deformation spaces of hyperbolic 3-manifolds.  Indeed,
Theorem~\ref{theorem:smooth:drill} represents the main technical tool
in the recent topological tameness theorems of the authors' with
R. Evans and J. Souto for algebraic limits of geometrically finite
manifolds, and the consequent reduction of Ahlfors' measure conjecture
to Conjecture~\ref{conjecture:density} (see
\cite{Ahlfors:finiteness}, \cite{BBES:tame}).

\bold{Grafting and geometric finiteness.}  Initially, our argument
mirrors that of \cite{Bromberg:bers}, in which a singly degenerate $M$
with arbitrarily short geodesics is first shown to be approximated by
geometrically finite cone-manifolds.  

The grafting construction of
\cite{Bromberg:bers} produces cone-manifolds that approximate a doubly
degenerate manifold as well, but the proof that these cone-manifolds
are geometrically finite is entirely different in this case.  Here, we
replace considerations of projective structures on surfaces with
notions of convex hulls and geometric finiteness for variable
(pinched) negative curvature developed by B. Bowditch and M. Anderson,
after applying a theorem of Gromov and Thurston to perturb the
relevant cone-metrics to smooth metrics of negative curvature.

\bold{Bounded geometry and arbitrarily short geodesics.}  After \cite{Bromberg:bers}, our central
challenge here is to address the possibility that $M$ is {\em doubly
degenerate}, namely, the case for which $M \cong S \times \reals$ and
the convex core is all of $M$.  In this case, $M$ has two {\em
degenerate ends}: each end has an exiting sequence of closed geodesics
that are homotopic to simple curves on $S$.  Our analysis turns on
whether such geodesics can be taken to be arbitrarily short.

When each end of $M$ has such a family of arbitrarily short geodesics,
a streamlined argument exists that avoids certain technical tools
developed here.  We refer the reader to
\cite[Sec. 3]{Brock:Bromberg:warwick} for a 
discussion of the argument, which is more directly analogous to that
of \cite{Bromberg:bers}.   We remark that in particular no application
of Thurston's double limit theorem is required; the convergence of the
relevant approximates follows directly from the cone-deformation theory.

When $M$ is assumed to have {\em bounded geometry} ($M$ has a global
lower bound to its injectivity radius) and $M$ is homotopy equivalent
to a surface, Minsky's {\em ending lamination theorem} for bounded
geometry implies Theorem~\ref{theorem:main} (see 
\cite[Cor. 2]{Minsky:bounded}).  The theorem guarantees that any such
$M$ is completely determined by its {\em end-invariants}, asymptotic
data associated to the ends of $M$.  An application of Thurston's
double limit theorem (\cite{Thurston:book:GTTM},
cf. \cite{Ohshika:ends}) and continuity of the length function for
laminations (see \cite{Brock:length}) allows one to realize the
end-invariants of $M$ as those of a limit $N$ of geometrically finite
manifolds $Q_n$.  Minsky's theorem \cite[Cor. 1]{Minsky:bounded} then
implies $N$ is isometric to $M$, and thus $\{Q_n\}$ converges to $M$.

A persistently difficult case has been that of $M$ with {\em
mixed-type}.  In this case, one end of $M$ has bounded geometry, the
other arbitrarily short geodesics.  For manifolds of mixed-type, the full
strength of our techniques is required to isolate the geometry of the
ends from one another.  Rather than breaking the argument into the
above cases, however, we have presented a unified treatment that
handles all cases simultaneously.

\bold{Scheme of the proof.}
As a guide to the reader, we briefly describe the scheme of the proof of
Theorem~\ref{theorem:main}.

\ital{I. Reduction to surface groups.}
The essential difficulties arise in
the search for geometrically finite approximates to a hyperbolic
3-manifold $M$ with the homotopy type of a surface $S$.  Within this 
category, it is the {\em doubly degenerate} manifolds that remain
after \cite{Bromberg:bers}.  Each such manifold has a positive and
a negative degenerate end, given a choice of orientation.

\ital{II. Realizing ends on Bers boundary.}  We first seek to
{\em realize} the geometry of each end of $M$ as that of an end of a
singly degenerate limit of quasi-Fuchsian manifolds
$\{Q(X,Y_n)\}$ or $\{Q(X_n,Y)\}$: given the positive end $E$ of $M$,
say, we seek a limit $Q = \lim Q(X,Y_n)$ so that $E$ admits a marking
and orientation preserving bi-Lipschitz diffeomorphism to an end of
$Q$.  We prove such limits can always be found
(Theorem~\ref{theorem:ends:realized}) by considering the 
bounded geometry case and the case when $E$ has arbitrarily short
geodesics separately.

\ital{III. Bounded geometry.} If the end $E$ has a lower
bound to its injectivity radius, we employ techniques of Minsky
to show its {\em end-invariant} $\nu(E)$ has {\em bounded type}: any
incompressible end of a hyperbolic 3-manifold with end-invariant
$\nu(E)$ has a lower bound to {\em its} injectivity radius, whether or
not the bound holds globally.  After producing
a limit $Q$ in a Bers boundary with end-invariant $\nu(E)$, an
application of Minsky's bounded geometry theory shows that $Q$ realizes
$E$ in the above sense.

\ital{IV. Arbitrarily short geodesics.} If the end $E$ has
arbitrarily short geodesics, a 
simultaneous grafting procedure produces a hyperbolic cone-manifold
with two cone-singularities with angle $4 \pi$.  Generalizing tameness
results for variable negative curvature, we show the simultaneous
grafting is {\em geometrically finite}: its convex core is compact.
Applying the drilling theorem (Theorem~\ref{theorem:bilip}) we
deform the metric back to a smooth structure {\em rel} the
conformal boundary with bounded distortion of the metric structure
outside tubes around the cone-singularities.  Successive
simultaneous graftings give quasi-Fuchsian manifolds
limiting to a manifold $Q$ that realizes $E$.

\ital{V. Asymptotic isolation.}  We then 
prove an asymptotic isolation theorem
(Theorem~\ref{theorem:isolation}) which again uses the drilling theorem to
show that any cusp-free doubly degenerate limit $M$ of quasi-Fuchsian
manifolds $Q(X_n,Y_n)$ has positive and negative ends $E^+$ and $E^-$
so that $E^+$ depends only on $\{Y_n\}$ and $E^-$ depends only on
$\{X_n\}$ up to bi-Lipschitz diffeomorphism.  

\ital{VI. Conclusion.}  The proof is concluded by realizing the
positive end $E^+$ of $M$ by the limit of $\{Q(X,Y_n)\}$ and the
negative end $E^-$ of $M$ by the limit of $\{ Q(X_n,Y)\}$ where
$\{X_n\}$ and $\{Y_n\}$ are determined by
Theorem~\ref{theorem:ends:realized}.  Thurston's double limit
theorem implies that $Q(X_n,Y_n)$ converges up to subsequence to a
limit $M'$, and thus Theorem~\ref{theorem:isolation} implies the ends
of $M'$ admit marking-preserving bi-Lipschitz 
diffeomorphisms to the ends of $M$.  By an application of Sullivan's
rigidity theorem, we have $Q(X_n,Y_n) \to  M$.

\bigskip

We conclude with two remarks.

\ital{Generalizations.}  The hypotheses of the theorem can
be weakened with only technical changes to the argument.  The clearly
essential hypothesis is that $M$ be {\em tame} which is guaranteed in
our setting by the assumption that $M$ have incompressible ends (by
Bonahon's theorem
\cite{Bonahon:tame}).  
In the setting of tame manifolds with compressible ends the principal
obstruction to carrying out our argument lies in the need for {\em
unknotted} short geodesics, guaranteed in the incompressible setting
by a result of J.P. Otal (see
\cite{Otal:unknotted2} and Thm.~\ref{theorem:otal}).  
We expect this to be a surmountable difficulty and will take up the
issue in a future paper

The assumption that $M$ have no parabolics is required only by our use
of Minsky's ending lamination theorem for bounded geometry
\cite{Minsky:ends} where hyperbolic manifolds in question are assumed
to have a global lower bound on their injectivity radii rather than
simply a lower bound to the length of the shortest geodesic.

A reworking Minsky's theorem to allow {\em
peripheral} parabolics represents the only obstacle to allowing
parabolics in our theorem.  While such a reworking is now
essentially straightforward after the techniques introduced in
\cite{Minsky:bounded}, we have chosen in a similar spirit to defer
these technicalities to a later paper in the interest of conveying the
main ideas.

\ital{Ending laminations.}  We also remark that recently announced
work of the first author with R. Canary and Y. Minsky \cite{BCM:elc}
has completed Minsky's program to prove Thurston's {\em ending
lamination conjecture} for hyperbolic 3-manifolds with incompressible
ends.  This result predicts (in particular) that each hyperbolic
3-manifold $M$ equipped with a cusp preserving homotopy equivalence
from a hyperbolic surface $S$ is determined up to isometry by its
parabolic locus and its end-invariants (see
\cite{Minsky:CKGI,BCM:elc}).

As in the bounded geometry case, Theorem~\ref{theorem:main} follows
from the ending lamination conjecture via an application of
\cite{Thurston:book:GTTM}, \cite{Ohshika:ends} and
\cite{Brock:length}, so the results of \cite{BCM:elc} will give an
alternate proof of our main theorem.  We point out that the techniques
employed here are independent of those of \cite{BCM:elc} and of a
different nature.  In particular, we expect the drilling theorem
(Theorem~\ref{theorem:bilip}) to have applications beyond the scope of
this paper, and we refer the reader to \cite{BBES:tame} for an initial
example of its application in a different context.

\bold{Acknowledgements.}  The authors are indebted to Dick Canary,
Craig Hodgson, Steve Kerckhoff, and Yair Minsky for their interest and
inspiration.

\section{Preliminaries}
\label{section:preliminaries}
A {\em Kleinian group} is a discrete, torsion free
subgroup of $\Isom^+(\half^3) = \Aut(\chat)$.  Each Kleinian group
$\Gamma$ determines a complete hyperbolic 3-manifold $M = \half^3 /
\Gamma$ as the quotient of $\half^3$ by $\Gamma$.
The manifold $M$ extends to its {\em Kleinian manifold} $N = (\half^3
\cup \Omega)/\Gamma$ by adjoining its {\em conformal boundary} $\bdry
M$, namely, the quotient by $\Gamma$ of the {\em domain of
discontinuity} $\Omega \subset \chat$ where $\Gamma$ acts properly
discontinuously.  (Unless explicitly stated, all Kleinian groups will
be assumed {\em non-elementary}).

The {\em convex core} of $M$, which we denote by $\core(M)$, is the smallest
convex subset of $M$ whose inclusion is a homotopy equivalence.  The
complete hyperbolic 3-manifold $M$ is {\em geometrically finite} if
$\core(M)$ has a finite volume unit neighborhood in $M$.

\bold{The thick-thin decomposition.}  The {\em injectivity radius} 
$\inj \colon M \to \reals^+$ measures the radius of the maximal
embedded metric open ball at each point of $M$.  
For $\epsilon >0$, we denote by $M^{< \epsilon}$ the $\epsilon$-thin
part where $\inj(x) <\epsilon$ and by $M^{\ge \epsilon}$ the {\em
$\epsilon$-thick part} $M \setminus M^{< \epsilon}$.
By the {\em Margulis
lemma} there is a universal constant $\varepsilon$ (depending only on
the dimension) so that each
component $T$ of the {\em thin part} $M^{< \varepsilon}$ where
$\inj(x) < \varepsilon$ has a standard type: either $T$ is an open
solid torus neighborhood of a short geodesic, or $T$ is the quotient
of an open horoball $B \subset \half^3$ by a $\zed$ or $\zed \dirsum
\zed$ parabolic group fixing $B$.  

\bold{Curves and surfaces.} Let $S$ be a closed topological surface of 
genus at least $2$.  We denote by $\eS$ the set of all isotopy classes
of essential simple closed curves on $S$.  The {\em geometric
intersection number} $$i \colon \eS \times \eS \to \zed^+$$ counts the
minimal number of intersections of representatives of curves in a pair
of isotopy classes $(\alpha,\beta) \in \eS \times \eS$.

The {\em Teichm\"uller space} $\Teich(S)$
parameterizes marked hyperbolic structures on $S$: pairs $(f, X)$
where $f \colon S \to X$ is a homeomorphism to a hyperbolic surface
$X$ modulo the equivalence relation that $(f,X) \sim (g,Y)$ when there 
is an isometry $\phi \colon X \to Y $ for which $\phi \compos f \simeq 
g$.   If we allow $S$ to have boundary, then $X$ is required to have
finite area and $f \colon \interior(S) \to X$ is a homeomorphism from
the interior of $S$ to $X$.

We topologize Teichm\"uller space by the {\em quasi-isometric distance}
$$d_{\rm qi}( (f,X), (g,Y) )$$ which is the log of the infimum over all
bi-Lipschitz diffeomorphisms $\phi \colon X \to Y$ homotopic to $g
\compos f^{-1}$ of the best bi-Lipschitz constant for $\phi$
(cf. \cite{Thurston:book:TDGT}).  Each $\alpha \in \eS$ has a unique
geodesic representative on any surface $(f,X) \in \Teich(S)$ by taking
the representative of the free-homotopy class of $f(\alpha)$ on $X$ of
shortest length.  

To interpolate between simple closed curves in $\eS$, Thurston
introduced the {\em measured geodesic lamiantions}, $\ml(S)$, which
may be obtained formally as the completion of the image of $\reals^+
\times \eS$ under the map
$\iota \colon \reals^+ \times \eS \to \reals^\eS$
defined by $\langle \iota(t,\alpha) \rangle_\beta = t
i(\alpha,\beta).$

On a given $(f,X)$ in Teichm\"uller space, a {\em geodesic lamination}
is a closed subset of $X$ given as a union of pairwise disjoint
geodesics on $S$.  The measured laminations $\ml(S)$ are then identified with
{\em measured geodesic laminations}, pairs $(\lambda,\mu)$ of a
geodesic lamination $\lambda$ and a {\em transverse measure}, an
association of a measure $\mu_\alpha$ to each arc $\alpha$ transverse
to $\lambda$ so that $\mu_\alpha$ is invariant under holonomy.
One obtains the {\em projective measured laminations} $\pl(S)$ as the
quotient $(\ml(S) \setminus \{0 \})/\reals^+$.  
(See
\cite{Thurston:book:GTTM}, \cite{FLP}, or \cite{Bonahon:currents} for
more about geodesic and measured laminations).

\bold{Surface groups.} 
By $H(S)$ we denote all marked hyperbolic 3-manifolds $(f \colon S \to
M)$: i.e. complete hyperbolic 3-manifolds $M$ equipped with homotopy
equivalences $f \colon S \to M$, modulo the equivalence relation $$(f
\colon S \to M) \sim (g \colon S \to N)$$ if there is an isometry
$\phi \colon M \to N$ for which $\phi \compos f \simeq g$.

Each $(f \colon S \to M)$ in $AH(S)$ determines a
representation $$f_* =  
\rho \colon \pi_1(S) \to \Isom^+(\half^3)$$ well defined up to
conjugacy in $\Isom^+(\half^3) = \PSL_2(\cx)$.  We topologize $H(S)$
by the compact open topology on the induced representations, up to
conjugacy.  Convergence in this sense is known as {\em algebraic
convergence}; we equip $H(S)$ with this {\em algebraic topology} to obtain
the space $AH(S)$, the {\em algebraic deformation space}.

The subset $QF(S) \subset AH(S)$ denotes the {\em quasi-Fuchsian
locus}, namely, manifolds $(f \colon S \to Q)$ so that $Q$ is
bi-Lipschitz diffeomorphic to the quotient of $\half^3$ by a Fuchsian
group.  Such a quasi-Fuchsian manifold $Q$ simultaneously uniformizes
a pair $(X,Y) \in \Teich(S) \times \Teich(S)$ as its conformal
boundary $\bdry Q$; in our convention $X$ compactifies the negative
end of $Q(X,Y) \cong S \times \reals$ and $Y$ compactifies the
positive end ($\chat$ is assumed oriented so that
 the resulting identification of $\wt S$ with $\wt Y
\subset \chat$ is orientation preserving and while the identification
of $\wt X$ with $\wt S$ is orientation reversing; by our convention
$Q(Y,Y)$ is a Fuchsian manifold).

Bers exhibited a
homeomorphism $$Q
\colon \Teich(S)
\times \Teich(S) \to QF(S)$$ that assigns to the pair $(X,Y)$ the
quasi-Fuchsian manifold $Q(X,Y)$ simultaneously uniformizing $X$ and
$Y$.  The manifold $Q(X,Y)$ naturally
inherits a homotopy equivalence $f \colon S \to Q(X,Y)$ from the
marking on either of its boundary components, so the simultaneous
uniformization is naturally an element of $AH(S)$.

One obtains a {\em Bers slice} of the quasi-Fuchsian space $QF(S)$ by
fixing one factor in the product structure; we denote by $$B_X = 
\{X\} \times \Teich(S) \subset QF(S)$$ 
the Bers slice of quasi-Fuchsian
manifolds with $X$ compactifying their negative ends.
As one may fix the conformal boundary compactifying either the
positive or negative end, we will employ the notation
$$B_X^+ = \{X\} \times \Teich(S) 
\ \ \ \ \text{and} \ \ \ \ 
B_Y^- = \Teich(S) \times \{Y\} $$
to distinguish the two types of slices.

If $g \colon M \to N$ is a bi-Lipschitz diffeomorphism between
Riemannian $n$-manifolds, its {\em bi-Lipschitz
constant} $L(g) \ge 1$ is the infimum over all $L$ for which
$$\frac{1}{L}  \le \frac{|D(v)|}{|v|} \le L$$
for all $v \in TM$.

Following McMullen, (see \cite[Sec. 3.1]{McMullen:book:RTM})  we
define the {\em quasi-isometric distance} on $AH(S)$ by 
$$d_{\rm qi}((f_1,M_1),(f_2,M_2)) = \inf \log L(g)$$
where the inf is taken over all orientation-preserving bi-Lipschitz
diffeomorphisms $g \colon 
M_1 \to M_2$ for which $g \compos f_1$ is 
homotopic to $f_2$.  If there is no such diffeomorphism in the
appropriate homotopy class, then we say $(f_1,M_1)$ and $(f_2,M_2) $
have infinite quasi-isometric distance. 
The quasi-isometric distance is lower semi-continuous on $AH(S) \times
AH(S)$ (\cite[Prop. 3.1]{McMullen:book:RTM}).

\bold{Geometric and strong convergence.}
Another common and related notion of convergence of hyperbolic
manifolds comes from the {\em Hausdorff topology}, which we now describe.

A hyperbolic 3-manifold determines a Kleinian group only up to
conjugation.  Equipping $M$ with a unit orthonormal frame $\omega$ at
a basepoint $p$ (a {\em base-frame}) eliminates this ambiguity via the
requirement that the covering projection $$ \pi \colon
(\half^3,\widetilde{\omega}) \to (\half^3,\widetilde{\omega})/
\Gamma = (M,\omega)$$
sends the standard frame $\widetilde{\omega}$ at the origin in
$\half^3$ to $\omega$.

The based hyperbolic 3-manifolds 
$(M_n,\omega_n) = (\half^3,\widetilde{\omega})/\Gamma_n$ converge {\em geometrically} to
a {\em geometric limit} $(N,\omega) = 
 (\half^3,\widetilde{\omega})/\Gamma_G$ 
if $\Gamma_n$ converges to
$\Gamma_G$ in the geometric topology: 
\begin{itemize}
\item For each $\gamma \in \Gamma_G$ there are $\gamma_n \in \Gamma_n$
with $\gamma_n \to \gamma$.
\item If elements $\gamma_{n_k}$ in a subsequence $\Gamma_{n_k}$
converge to $\gamma$, then  $\gamma$ lies in $\Gamma_G$.
\end{itemize}

Geometric convergence has an internal formulation: $(M_n,\omega_n)$
converges to $(N,\omega)$ if for each smoothly embedded compact
submanifold $K \subset N$ containing $\omega$, there are
diffeomorphisms $\phi_n \colon K \to (M_n,\omega_n)$ so that
$\phi_n(\omega) = \omega_n$ and so that $\phi_n$ converges to an
isometry on $K$ in the $C^\infty$ topology
(\cite{Benedetti:Petronio:book} \cite[Ch. 2]{McMullen:book:RTM}).

When $(f \colon S \to M)$ lies in $AH(S)$, a base-frame $\omega \in M$
determines a discrete faithful representation $f_* \colon \pi_1(S) \to
\Gamma$ where $(M,\omega) = (\half^3,\widetilde{\omega})/\Gamma$.
Denote by $AH_\omega(S)$ the marked {\em based} hyperbolic 3-manifolds
$(f\colon S \to (M,\omega))$; based hyperbolic 3-manifolds
$(M,\omega)$ together with homotopy equivalences $f \colon S \to
 (M,\omega)$ up to isometries that preserve marking and base-frame.

The space $AH_\omega(S)$ carries the topology of convergence on
generators of the induced representations $f_*$; the topology on
$AH(S)$ is simply the quotient topology under the natural base-frame
forgetting map $AH_\omega(S) \to AH(S)$.  As with $AH(S)$ we will
often assume an implicit marking and refer to $(M,\omega) \in AH_\omega(S)$.

Consideration of $AH_\omega(S)$ allows us to understand the relation
between algebraic and geometric convergence (see 
\cite[Sec. 2]{Brock:iter}):
\begin{theorem}
Given a sequence $\{(f_n \colon S \to M_n)\}$ with limit
$(f \colon S \to M)$ in $AH(S)$ there are convergent lifts $(f_n
\colon S \to (M_n,\omega_n))$ to $AH_\omega(S)$ so that after passing
to a subsequence $(M_n,\omega_n)$ converge geometrically to a
geometric limit $(N,\omega)$ covered by $M$ by a local isometry.
\end{theorem}

When this local isometry is actually an isometry, we say the
convergence is {\em strong}.
\begin{defn}
The sequence $M_n \to M$ in $AH(S)$ converges {\em
strongly} if there are lifts 
$(M_n,\omega_n) \to (M,\omega)$ to $AH_\omega(S)$ so that
$(M_n,\omega_n)$ also converges geometrically to $(M,\omega)$.
\end{defn}

\bold{Pleated surfaces.}  Given an $M \in AH(S)$ and a simple closed curve $\alpha \in
\eS$ representing a non-parabolic conjugacy class of $\pi_1(M)$,
we follow Bonahon's convention and denote by $\alpha^*$ the geodesic
representative of $\alpha$ in $M$.  To control how $\alpha^*$ can lie
in $M$, Thurston introduced the notion of a {\em pleated surface}.
\begin{defn}
A path isometry $g \colon X \to N$ from a hyperbolic surface $X$ to a
hyperbolic 3-manifold $N$ is a {\em pleated surface} if for 
each $x \in X$ there is a geodesic segment $\sigma$ through $x$ so
that $g$ maps $\sigma$ isometrically to $N$.
\end{defn}

When $M$ lies in $AH(S)$, a particularly useful class of pleated
surfaces arises from those that ``preserve marking'' in the following
sense: denote by $\ps(M)$ the set of all pairs $(g,X)$ where $X$
lies in $\Teich(S)$ and $g \colon X \to M$ is a pleated surface with
the property that $g \compos \phi \simeq f$ where $\phi$ is the
implicit marking on $X$ and $f$ is the implicit marking on $M$.

Given a lamination $\mu \in \ml(S)$, we say the pleated surface $(g,X)
\in \ps(M)$ {\em realizes} $\mu$ if each geodesic leaf $\ell$ in the support of
$\mu$ realized as a geodesic lamination on $X$ is mapped by $g$ by a
local isometry; alternatively, the lift $\wt g \colon \wt X \to
\half^3$ sends each lifted leaf $\wt \ell$ of $\wt \mu$ to a complete
geodesic in $\half^3$.

The following {\em bounded diameter theorem} for pleated surfaces is
instrumental in Thurston's studies of geometrically tame hyperbolic
3-manifolds (see \cite[Sec. 8]{Thurston:book:GTTM} or alternate
versions in \cite{Bonahon:tame}, \cite{Canary:ends}).
\begin{theorem}
Each compact subset $K \subset M$ has an enlargement $K'$ so that if
$(g,X) \in \ps(M)$, and $g(X) \cap K \not = \nullset$ then $g(X)$ lies
entirely in $K'$.
\label{theorem:diameter}
\end{theorem}
 
\bold{Tame ends.}  Let $M$ be a complete hyperbolic 3-manifold with
finitely generated fundamental group.  By a theorem of P. Scott
\cite{Scott:core}, there 
is a compact submanifold $\calM \subset M$ whose inclusion is a
homotopy equivalence.  By convention, given a choice of compact core
$\calM$ for $M$, the {\em ends} of $M$ are the connected components of
the complement $M \setminus \calM$ of the interior of $\calM$.  Each
end $E$ is {\em cut off} by a boundary component $S \subset \bdry
\calM$.

The end $E$ is {\em tame} if it is homeomorphic to the product $S 
\times \reals^+$, and the manifold $M$ is {\em topologically tame} (or
simply {\em tame}) if it is the interior of a compact manifold with
boundary.  When $M$ is tame we can choose the compact core $\calM$
such that each end $E$ is tame. Manifestly, the end $E$ depends on a
choice of compact core, but as we will typically be interested in the
end $E$ only up to bi-Lipschitz diffeomorphism, we will assume such a
core to be chosen in advance and address any ambiguity as the need
arises. 

An end $E$ of $M$ is {\em geometrically finite} if it has finite
volume intersection with the convex core of $M$.  Otherwise it is {\em
geometrically infinite}.  By a theorem of Marden (see
\cite{Marden:kgs}), a geometrically finite end is tame.  A
geometrically infinite tame end is {\em 
simply degenerate}.  The manifold $M$ is geometrically finite if and
only if each of its ends is geometrically finite.

\bold{Tameness and Otal's theorem.}  One key element of our argument
in section~\ref{section:grafting}
involves the fact that any collection of sufficiently short closed
curves in $M \in AH(S)$ is {\em unknotted} and {\em unlinked}.

Given $M \in AH(S)$ the tameness theorem of Bonahon and Thurston
\cite{Bonahon:tame,Thurston:book:GTTM} guarantees the existence of a
product structure $F \colon S \times \reals \to M$.  Otal defines a
notion of `unknottedness' with respect to this product structure as
follows: a closed curve $\alpha \in M$ is {\em unknotted} if it is
isotopic in $M$ to a simple curve in a {\em level surface} $F(S \times
\{t\})$.  Likewise, a collection $\calC$ of closed curves in $M$ is {\em
unlinked} if there is an isotopy of the collection $\calC$ sending each
member $\alpha \in \calC$ to a distinct level surface $F(S \times
\{t_\alpha\})$.

\begin{theorem}[Otal]
Let $S$ be a closed surface, and let $(f \colon S \to M)$ lie in
$AH(S)$.  There is a constant $\ell_{\rm knot} >0$ depending only on
$S$ so that if $\calC$ is any collection of closed curves in $M$ for which
for which $\ell_M(\alpha^*)< \ell_{\rm knot}$ for each $\alpha \in
\calC$, then the collection $\calC$ is unlinked.
\label{theorem:otal}
\end{theorem}
See \cite{Otal:unknotted} and, in particular, \cite[Thm. B]{Otal:unknotted2}.

\bold{Cone-manfiolds.}
A key technical ingredient for our argument will be the notion of a
{\em 3-dimensional hyperbolic cone-manifold}.  Let $N$ be a compact
3-manifold with boundary and $\calC$ a collection of disjoint simple
closed curves. A {\em hyperbolic cone-metric} on $(N, \calC)$ is a
hyperbolic metric on the interior of $N \backslash \calC$ whose
completion is a singular metric on the interior of $N$. In a
neighborhood of a point in $\calC$ the metric will have the form
$$dr^2 + \sinh^2r d\theta^2 + \cosh^2 r dz^2$$ 
where $\theta$ is measured the cone-angle $\alpha$. The singular locus
will be identified with the $z$-axis and will be totally
geodesic. Note that the cone angle will be constant along each
component of the cone singularity.

\section{Geometric finiteness in negative curvature}
\label{section:geom:finite}
In this section, define the notion of geometric finiteness
for 3-dimensional hyperbolic cone-manifolds we will use and show its
equivalence to precompactness of the set of closed geodesics in the
cusp-free setting.  We then go on to employ the work of Bonahon and
Canary \cite{Bonahon:tame} \cite{Canary:ends} to show the existence of
simple closed geodesics exiting any end of $M$ that is not
geometrically finite.

\bold{Geometric finiteness for cone-manifolds.}  When the convex core
of the complete hyperbolic 3-manifold $M$ has a finite-volume unit
neighborhood, the only obstruction to the compactness of the convex
core is the presence of cusps in $M$.  In the cusped case, a slightly
different definition is required.  For our discussion, we consider
only cusps that arise from rank-two abelian subgroups of the
fundamental group, or {\em rank-two cusps}.
\begin{defn}
A 3-dimensional hyperbolic cone-manifold  $M$ is {\em  geometrically
finite without rank-one cusps} if $M$ has a compact core bounded by
convex surfaces and tori. 
\end{defn}
In the sequel, all hyperbolic cone-manifolds we will consider will be
free of rank-one cusps.  As such, we simply refer to geometrically
finite manifolds without rank-one cusps as
{\em geometrically finite}.

Given a compact core $\calM$ for such an $M$, the geometric
finiteness of $M$ is usefully rephrased as a condition on the
ends of $M$ (again, we refer to components of $M \setminus \calM$ as
the {\em ends} of $M$; they are neighborhoods of the topological ends
of $M$).   

\begin{defn}
An end $E$ of a 3-dimensional hyperbolic cone-manifold $M$
is {\em geometrically finite} if its intersection with the
convex core of $M$ has finite volume.
\end{defn}
An end that is cut off by a torus will be a rank two cusp and will be
entirely contained in the convex core. Since we are assuming $M$ does
not have rank-one cusps, each end of a geometrically finite manifold
cut off by a higher genus surface will intersect the convex core in a
compact set.

Then one may easily verify the following proposition.
\begin{prop}
The 3-dimensional hyperbolic cone-manifold $M$ is geometrically finite 
if and only if each end of $M$ is geometrically finite.
\end{prop}

\bold{Geometrically infinite ends.}
An end $E$ of $M$ that is not geometrically finite is {\em
geometrically infinite} or {\em degenerate}. 
\begin{defn}
Let $E$ be a geometrically infinite end of a 3-dimensional hyperbolic
cone-manifold $M$, cut off by a surface $S$.  Then $E$ is {\em simply
degenerate} if for any compact subset $K \subset E$ there is a simple
curve $\alpha$ on $S$ whose geodesic representative lies in $E \setminus K$.
\end{defn}
In the smooth hyperbolic setting, a synonym for a simply degenerate
end is a {\em geometrically tame} end; we use the same terminology
here.  The cusp-free hyperbolic cone-manifold $M$ is {\em
geometrically tame} if all its ends are geometrically finite or
geometrically tame.

Thurston and Bonahon proved that a geometrically tame manifold $M$ is
{\em topologically tame}, namely, $M$ is homeomorphic to the interior
of a compact 3-manifold.  Generalizing Bonahon's work, Canary proved
the converse:
\begin{theorem}[Canary]
Let $M$ be a topologically tame complete hyperbolic 
3-manifold.  Then $M$ is geometrically tame.
\end{theorem}

\bold{Geometric finiteness in variable negative curvature.}
Brian Bowditch has given a detailed analysis of how various notions of
geometric finiteness for complete hyperbolic 3-manifolds and their
equivalences generalize to the case of {\em pinched negative
curvature} namely, 3-manifolds with complete Riemannian metrics with
all sectional curvatures in the interval $[-a^2,-b^2]$, where $0 < b <
a$.

Such a manifold is the quotient of a {\em pinched Hadamard manifold}
$X$, a simply connected manifold with sectional curvatures pinched
between $-a^2$ and $-b^2$, by a discrete subgroup $\Gamma$ of its
orientation preserving isometries $\Isom^+ X$.  For our purposes, we
assume $X$ has dimension $3$.  The action of $\Gamma$ on $X$ has much
in common with actions of Kleinian groups on $\half^3$.  In
particular, $X$ has a natural {\em ideal sphere} $X_I$ or {\em sphere
at infinity}, which may be identified with equivalence classes of
infinite geodesic rays in $X$ where rays are equivalent if they are
asymptotic.

As in the hyperbolic setting, the action of $\Gamma$ on $X_I$ is
partitioned into its {\em limit set} $\Lambda$ where the orbit of a (and hence
any) point in $X$ accumulates on $X_I$ and its {\em domain of
discontinuity} $\Omega  =  X_I \setminus \Lambda$.

The {\em convex core} of the quotient manifold $M = X/\Gamma$ of
pinched negative 
curvature is the quotient $\hull(\Lambda) /\Gamma$ of the convex hull
in $X$ of the limit set $\Lambda$ by the action of
$\Gamma$.  Then following \cite{Bowditch:finiteness} we make the
following definition.
\begin{defn}
The manifold $M = X/\Gamma$ of pinched negative curvature is {\em
geometrically finite} if the radius-$1$ neighborhood of the 
convex core has finite volume.
\end{defn}
In the case when the manifold $M$ has no cusps (the group $\Gamma$ is
free of parabolic elements; see \cite[Sec. 2]{Bowditch:finiteness}),
this notion is equivalent to the compactness of the convex core.

Indeed, it suffices to consider the quotient of the {\em join} of the
limit set $\join(\Lambda)$: the collection of all geodesics in $X$
joining pairs of points in $\Lambda$.
We apply the following theorem of Bowditch
\cite{Bowditch:hulls} which follows from work of M. Anderson
\cite{Anderson:dirichlet}.
\begin{theorem}[Bowditch]
Let $M$ be a Riemannian manifold of pinched negative curvature.  Then
there is a $\sigma>0$ depending only on the pinching constants so that
$$\hull(\Lambda) \subset \calN_\sigma(\join(\Lambda)).$$
\label{theorem:Bowditch:hull}
\end{theorem}
(Cf. \cite[Sec. 5.3]{Bowditch:finiteness}).

In the complete smooth hyperbolic setting, the density of the fixed
points of hyperbolic isometries in $\Lambda \times \Lambda$ gives
another characterization of geometric finiteness: the complete
cusp-free hyperbolic 3-manifold is geometrically finite if and only if
the closure of the set of closed geodesics in $M$ is compact.

A lacuna in the various existing discussions of how features of the
complete hyperbolic setting generalize to the pinched negative
curvature setting is the following equivalence, which will allow us to
improve the results of Canary \cite{Canary:ends}.
\begin{lem}
Let $M$ be a 3-dimensional manifold of pinched negative curvature and
no cusps.  Then $M$ is geometrically finite if and only if the closure
of the set of closed geodesics in $M$ is compact.
\label{lemma:pinched:gf}
\end{lem}

\bold{Proof:}  Let $M = X /\Gamma$, where $X$ is a 3-dimensional
pinched Hadamard manifold and $\Gamma$ is a discrete subgroup of
$\Isom^+ X$.

Since fixed points of hyperbolic isometries are again dense in $\Lambda \times
\Lambda$, it follows that lifts of closed geodesics to $X$ are dense in
$\join(\Lambda)$.  Applying Theorem~\ref{theorem:Bowditch:hull}, we
have $$\hull(\Lambda)
\subset\calN_{\sigma}(\join(\Lambda))$$ where $\sigma$ depends
only on the pinching constants for $M$.
But if the closure of the set of closed geodesics in $M$ is compact
then the quotient $\closure{\join(\Lambda)}/\Gamma$ is compact.
It
follows that the convex core 
$$\core(M) \subset \calN_{\sigma}(\closure{\join(\Lambda)}/\Gamma)$$
is compact.
Thus $M$ is a
geometrically finite manifold of pinched negative curvature.  

Conversely, since all closed geodesics in $M$ lie in $\core(M)$, 
the closure of the set of closed geodesics in $M$ is compact whenever
$\core(M)$ is compact. 
\qed

\begin{cor}
Let $M$ be a 3-dimensional hyperbolic cone-manifold with no cusps
so that for every cone-singularity $c$ the cone-angle at $c$ is
greater than $2 \pi$.  Then $M$ is geometrically finite if and only if
the closure of the set of all closed geodesics in $M$ is compact.
\end{cor}

\bold{Proof:}
By a standard argument (see \cite{Gromov:Thurston:pinching}) the
assumption on the cone-angles implies that the singular hyperbolic
metric on $M$ may be perturbed to give a negatively curved metric on
$M$ that is hyperbolic away from a tubular neighborhood of the
cone-locus.  The result is a Riemannian manifold of pinched negative 
curvature $\hat{M}$.

The smoothing $\hat{M}$ is a new metric on $M$ and in this new metric
each closed geodesic is a uniformly bounded distance from its geodesic
representative in $M$.
It follows that the closure of the set of closed geodesics in $M$ is
compact if and only if the closure of the set of closed geodesics in
$\hat{M}$ is compact.  

If the union of all closed geodesics in $\hat M$ is precompact, then
$\hat M$ is geometrically finite, by Lemma~\ref{lemma:pinched:gf}.  It
follows that the convex core for $\hat M$ is compact (since $\hat M$
has no cusps).  For $R>0$ sufficiently large, the radius-$R$
neighborhood of the convex core of $\hat M$ gives a compact core
$\calM$ for $\hat M$ bounded by convex surfaces that miss the
neighborhoods where the metrics on $M$ and $\hat M$ differ.  Since
convexity is a local property for embedded surfaces, it follows that
the surfaces $\bdry \calM$ are convex in $M$ as well.

We conclude that if the closed geodesics are precompact in $M$ then
$M$ has a compact core bounded by convex surfaces, so $M$ is
geometrically finite.  The converse is immediate.
\qed

We now prove the appropriate generalization of Canary's theorem in the 
context of hyperbolic cone-manifolds.  As in the smooth case, we say a
3-dimensional hyperbolic cone-manifold is {\em topologically tame} if it is
homeomorphic to the interior of a compact 3-manifold.
\begin{theorem}
Suppose $M$ is a topologically tame 3-dimensional hyperbolic
cone-manifold.  Assume that the cone-angle at each cone-singularity is
at least $2 \pi$.  Then $M$ is geometrically tame: each end $E$ of $M$
is either geometrically finite or simply degenerate.
\end{theorem}

\bold{Proof:}
We will apply Canary's generalizations of Bonahon's arguments
(\cite{Bonahon:tame,Canary:ends})  to the
variable negative curvature setting, after
taking care to ensure that in a geometrically infinite end $E$ one can
always find closed geodesics (not necessarily simple) exiting $E$.
This will follow from an argument of Bonahon \cite{Bonahon:tame}.

\ital{Closed geodesics lie in every neighborhood of $E$.}
As above, we let $\hat M$ be a smoothing of $M$ to a manifold of
pinched negative curvature, modifying the metric in a close
neighborhood of the singular locus.  Since neighborhoods of the ends
are unchanged by this smoothing, it suffices to prove the theorem for
$\hat M$.

Let $E$ be a geometrically infinite end of $\hat{M}$ cut off by a
surface $S_0$, and let $K$ be a compact submanifold of $E$ so that 
$\bdry K = S_0 \disjunion S$ where $S$ is a smooth surface in $E$.  
We claim that there exists closed curves on $S_0$ whose geodesic
representatives eventually lie outside of $K$.  

Assuming otherwise, consider a sequence $\gamma_n$ of closed curves on 
$S_0$ whose geodesic representatives all intersect $K$, but whose
closure is not compact.  Then for any compact subset $K'$ in $E$ there 
is a $\gamma_n$ for which some arc of intersection of $\gamma_n^*$
with $E \setminus K$ exits $K'$.  Consider a family of such arcs
$a_n$ of intersection of $\gamma_n^*$ with $E \setminus K$ that eventually
intersect the complement of each compact subset of $E$.  The 
unordered pairs of endpoints $(x_n,y_n)$ of $a_n$ range in the compact
set $S \times S/\{(x,y) \sim (y,x) \}$.  

Pass to a subsequence of $\{(x_n,y_n)\}$ converging to
$(x_\infty,y_\infty)$.  Given $\epsilon>0$, we may choose an arc $a
\in \{a_n\}$ whose endpoints are within $\epsilon$ of
$(x_\infty,y_\infty)$.  We may choose another arc $a' \in \{a_n\}$
whose endpoints are also within $\epsilon$ of $(x_\infty,y_\infty)$,
and whose length is at least four times that of $a$.  Let $b$ be the
closed loop obtained from joining the endpoints of $a$ with those of
$a'$ by short arcs on $S$ of length less than $2\epsilon$.

Since the length of $a'$ in $\hat{M}$ is more than twice as long as
the combined length of the other three edges of $b$ (taking $\epsilon$
small), the lift of $b$ to the universal cover is a uniform
quasi-geodesic (see \cite{Bridson:Haefliger:npc}).  Thus, $b$ is in
particular homotopically non-trivial, and its geodesic representative
$b^*$ lies a uniformly bounded Hausdorff distance $d$ from $b$.

The loop $b$ is homotopic to a closed curve on $S$, and its geodesic
representative lies in the uniform neighborhood $\calN_d(E \setminus K)$ of the
complement of $K$ in $E$.  Since $K$ is arbitrary and $d$ does not
depend on $E$ or $K$, there are closed loops $\alpha_n$ so that given
any compact subset $K$ of $E$, there is an $n$ for which the geodesic
representative $\alpha_n^*$ of $\alpha_n$ lies in $E \setminus K$.

\ital{Topologically tame implies geometrically tame.}
In \cite{Canary:ends} a generalization of Bonahon's tameness theorem
\cite{Bonahon:tame} is applied in the context of branched covers of
hyperbolic 3-manifolds.  
After smoothing the branching locus to obtain a manifold with pinched
negative curvature that is hyperbolic outside of a compact set, Canary 
discusses the appropriate generalization to the main theorem of
\cite{Bonahon:tame} in this context.
In \cite[Sec. 4]{Canary:ends}, however, a {\em geometrically infinite end $E$} cut
off by $S$ is defined to be an end for which there are closed loops $\alpha_n
\subset S$ whose geodesic representatives eventually lie outside of
every compact subset of $E$.  
The above shows that if an end $E$ is not geometrically finite in our
sense, then it is geometrically infinite in the sense of Canary
\cite[Sec. 4]{Canary:ends}. 

Applying the tameness theorem of \cite[Thm. 4.1]{Canary:ends}, then,
if $M$ is a tame 3-dimensional hyperbolic cone-manifold, then all of
its ends are either geometrically finite or simply degenerate.
\qed

\section{Bounded geometry}
\label{section:bounded:geometry}
A central dichotomy in the study of ends of hyperbolic 3-manifolds
lies in the distinction between hyperbolic manifolds with {\em bounded
geometry} and those with {\em arbitrarily short geodesics}.  

A recent theorem of Y. Minsky shows that whether a manifold $M \in
AH(S)$ has bounded geometry is predicted by a comparison of its
end-invariants, a collection of geodesic laminations and hyperbolic
surfaces associated to the ends of $M$.

In this section we adapt Minsky's techniques to produce a version of
these criteria which can be 
applied end-by-end: we show that whether or not a simply degenerate end $E$
has bounded geometry depends only on {\em its} ending lamination $\nu(E)$
and not on the remaining ends.

\bold{End-invariants.}  When an end $E$ of $M \in AH(S)$ is
geometrically finite it admits a foliation by surfaces whose geometry
is exponentially expanding, but whose conformal structures converge to
that of a component, say $X$, of the conformal boundary of $M$.  With
the induced marking from $f$, $X$ determines a point in Teichm\"uller
space, and the asymptotic geometry of the end $E$ is determined by
this marked Riemann surface.  We say $X$ is the {\em end-invariant} of
the geometrically finite end $E$ (see, e.g.,
\cite{Epstein:Marden:convex} \cite{Minsky:ends}).   

A simply degenerate end of $M$ also has a well defined 
end-invariant.
\begin{defn}
Let $E$ be a simply degenerate end of $M$ cut off by a surface
$S$.  Let $\alpha_n$ be a sequence of simple closed curves on $S$
whose geodesic representatives $\alpha_n^*$ leave every compact subset 
of $E$.  Then the support $|[\nu]|$ of any limit $[\nu] \in \pl(S)$ of 
$\alpha_n$ is the {\em ending lamination} of $E$.
\end{defn}
By a theorem of Thurston, any two limits $[\nu]$ and $[\nu']$ in
$\pl(S)$ satisfy 
$$|\nu | = |\nu'|$$ so $\nu(E)$ is well defined.
We call the ending lamination $\nu(E)$ the {\em end-invariant}
for the degenerate end $E$.

For each $M$ in $AH(S)$ with no cusps, we will denote by $\nu^-$
and $\nu^+$ the end-invariants of the positive and negative ends $E^-$
and $E^+$ of $M$.  
\smallskip

\bold{Curve complexes and projections.}
In \cite{Harvey:CC}, W. Harvey organized the simple closed curves on
$S$ into a complex in order to develop a better understanding of the
action of the mapping class group.  Recently, (see
\cite{Minsky:bounded,Minsky:KGCC,Brock:wp})
his complex has become a fundamental object in the
study of 3-dimensional hyperbolic manifolds.

The {\em complex of curves} $\calC(S)$ is obtained by associated a
vertex to each element of $\eS$ and stipulating that $k+1$ vertices
determine a $k$-simplex if the corresponding curves can be realized
disjointly on $S$.  Except for some sporadic low genus cases, the same
definition works for non-annular surfaces with boundary (provided
$\eS$ is taken to represent the isotopy classes of {\em non-peripheral}
essential simple closed curves on $S$) and a similar {\em arc-complex}
can be defined for consideration of the annulus.  A remarkable theorem
of H. Masur and Y. Minsky establishes that the natural distance on the
$\calC(S)$ obtained by making each $k$-simplex a standard Euclidean
simplex turns $\calC(S)$ into a $\delta$-hyperbolic metric space (see
\cite{Masur:Minsky:CCI} for more details).

When $Y \subset S$ is a proper essential subsurface of $S$, $\calC(Y)$
is naturally a subcomplex of $\calC(S)$.  Masur and Minsky define a
projection map $\pi_Y \colon \calC(S) \to \calP(\calC(Y))$ from
$\calC(S)$ to the set of subsets of $\calC(Y)$ by associating to each
$\alpha \in \calC(S)$ the arcs of essential intersection of $\alpha$
with $Y$, surgered along the boundary of $Y$ to obtain simple closed
curves in $Y$.  The possible surgeries can produce curves in $Y$ that
intersect, but given any simplex $\sigma \in \calC(S)$ the total
diameter of $\pi_Y(\sigma)$ in $\calC(Y)$ is at most $2$ (see
\cite[Sec. 2, Lem. 2.3]{Masur:Minsky:CCII} for more details).

The {\em projection distance} $d_Y(\alpha,\beta)$ 
measures the distance from $\alpha$ to $\beta$
relative to the subsurface $Y$: $$d_Y(\alpha,\beta) =
\{\diam_{\calC(Y)}(\pi_Y(\alpha) \cup \pi_Y(\beta))\}.$$
Note that by the above, the projection $\pi_Y$ is 2-Lipschitz, i.e. we
have
$$d_{\calC(S)}(\alpha,\beta) \le 2 d_Y(\alpha,\beta)$$
for any pair of vertices $\alpha$ and $\beta$ in $\calC(S)$.

By a result of E. Klarreich \cite{Klarreich:boundary}, the Gromov
boundary of $\calC(S)$ is in bijection with the
possible ending laminations for a cusp-free simply degenerate end of
$M \in AH(S)$.  We denote this collection of geodesic laminations by
$\el(S)$.  Given such an ending lamination $\nu$, the projection
$\pi_Y(\nu)$ can be defined just as for $\alpha \in
\calC(S)$, and $\pi_Y(\nu)$ is the limiting value of $\pi_Y(\alpha_i)$
where $\alpha_i$ converges to $\nu \in \bdry \calC(S)$. 

If $Z \in \Teich(S)$ is a conformal boundary component of $M \in
AH(S)$, there is a uniform upper bound to the length of the shortest
geodesic on $Z$.  Although the shortest geodesic may not be unique,
the set $\short(Z)$ of shortest geodesics on $Z$ determines a set of
uniformly bounded diameter in $\calC(S)$.  Thus, given end-invariants
$\nu^-$ and $\nu^+$ for a cusp free $M\in AH(S)$, we can compare the
end-invariants in the surface $Y$ by the quantity $$d_Y(\nu^-,\nu^+)$$
where if $\nu^- = Z \in \Teich(S)$ we replace $\nu^-$ with
$\short(Z)$.

Using such comparisons, the main results of
\cite{Minsky:KGCC} and
\cite{Minsky:bounded} give necessary and sufficient conditions for the
length of the shortest closed geodesic in $M$ to have a lower bound
$\ell_0>0$.
\begin{theorem}[Minsky]
Let $M\in AH(S)$ have no cusps and end-invariants $(\nu^-,
\nu^+)$.  Then $M$ has bounded geometry if and only if the supremum
$$\sup_{Y \subset S} d_Y(\nu^-, \nu^+)$$
over all proper essential subsurfaces $Y \subset S$ is bounded above.
\label{Minsky:bounded}
\end{theorem}

We deduce the following corollary.
\begin{cor}
Let the doubly degenerate manifold $M \in AH(S)$ have no cusps.  
If the positive end $E^+$ of $M$ has 
bounded geometry, any degenerate manifold $Q$ in the Bers slice
$B_Y$ with ending lamination $\nu(E^+)$ has bounded geometry.
\label{corollary:bounded:type}
\end{cor}

\bold{Proof:}
Assume otherwise.  Then by Minsky's theorem, there exists a family of
essential subsurfaces $Y_j \subset S$ so that 
$$d_{Y_j}(Y, \nu^+) \to \infty$$ as $j$ tends to $\infty$.  Choosing
$\alpha_j \subset \bdry Y_j$, we have 
$$\ell_Q(\alpha_j) \to 0$$ by
\cite[Thm. B]{Minsky:KGCC}.
Since the geodesic representatives
$\alpha_i^*$  in $Q$ exit the end of $Q$, any limit $[\nu]$ of
$\alpha_i$ in
$\pl(S)$ has intersection number zero with $\nu^+$ by the exponential
decay of the intersection number 
(see \cite[Ch. 9]{Thurston:book:GTTM}, \cite[Prop. 3.4]{Bonahon:tame}).  

We claim that the projection sequence
$$\{ d_{Y_j}(\nu^-,\nu^+) \}$$ is also unbounded.

Consider the distances
$$d_{Y_j}(\nu^-, Y).$$ 
Then either
$d_{Y_j}(\nu^-, Y)$ remains bounded or we may pass to a subsequence
so that $d_{Y_j}(\nu^-, Y) \to \infty$.

In the first case we have by the triangle inequality,
$$d_{Y_j}(Y,\nu^+) \le d_{Y_j}(Y, \nu^-) +  d_{Y_j}(\nu^-,\nu^+),$$
in particular, $d_{Y_j}(\nu^-,\nu^+)$ is unbounded.
By the main theorem of \cite{Minsky:KGCC}, it
follows that the simple closed  curves $\alpha_i$ satisfy
$\ell_M(\alpha_i) \to 0$.  

Thus, the geodesic representatives
$\alpha_i^{**}$  of $\alpha_i$ in $M$ must exit the end $E^-$ of $M$,
since their lengths have zero infimum.  Again applying
\cite[Ch. 9]{Thurston:book:GTTM}, \cite[Prop. 3.4]{Bonahon:tame}, we have
$[\nu]$ has intersection number zero with $\nu^-$.  It follows that
$\nu^- = \nu^+$, a contradiction (the ending laminations of a
cusp-free doubly degenerate manifold $M$ in $AH(S)$ must be distinct.
See \cite[Sec. 5]{Bonahon:tame} \cite[Ch. 9]{Thurston:book:GTTM}).
Thus, $Q$ has bounded geometry in this case.

If, on the other hand $d_{Y_j}(\nu^-, Y) \to \infty$, we consider 
a limit $Q'$ of quasi-Fuchsian manifolds in the Bers slice $B_Y$ with
ending lamination $\nu^-$.  Then by
\cite{Minsky:KGCC}
the curves $\alpha_j$ again have the property that
$\ell_{Q'}(\alpha_j) \to 0$, so the geodesic representatives
of $\alpha_j$ in $Q'$ must exit the end of $Q'$.  We again arrive at
the contradiction $\nu^- =\nu^+$, so we may conclude again that
$Q$ has bounded geometry.
\qed

The argument motivates the following definition.
\begin{defn}
A lamination $\nu \in \calE\calL(S)$ has {\em bounded type} if for any
$\alpha \in \calC(S)$, 
$\{d_{Y_j}(\alpha,\nu)\}$ is bounded over all essential subsurfaces $Y_j
\subset S$.
\end{defn}

\bold{Remark:}  The projection distances $d_{Y_j}(\alpha,\nu)$ are
reminiscent of the continued fraction expansion of an irrational
number.  In the case when $S$ is a punctured torus, this analogy is
literal in the sense that simple closed curves on $S$ are encoded by
their rational slopes, and measured laminations (up to scale) are
naturally the completion of the simple closed curves (see
\cite{Minsky:torus}).  In the punctured torus setting, bounded type
laminations are encoded by {\em bounded type} irrationals, namely,
irrationals with uniformly bounded continued fraction expansion.

\begin{theorem}
Let $E$ be a geometrically infinite tame end of a cusp-free hyperbolic
3-manifold $M$, and assume there is a lower bound to the injectivity
radius on $E$.  
Then there is a compact set $K \subset E$ and a
manifold $Q_\infty$ in the Bers boundary $\bdry B_Y$  so that the
subset $E \setminus K$ is bi-Lipschitz diffeomorphic to a the complement
$E_\infty \setminus K_\infty$ of a compact subset $K_\infty \subset Q_\infty$.
\label{theorem:bounded:geometry}
\end{theorem}

\bold{Proof:}  Let $\nu = \nu(E)$ be the ending lamination
for $E$.  Since the injectivity radius of $E$ is bounded below, it
follows that $\nu$ has bounded type.  By an application of the
continuity of the length function for laminations on $AH(S)$,
\cite[Thm. 1.3]{Brock:boundaries} there exists some $Q_\infty \in
\bdry B_Y$ so that $\nu(Q_\infty) =
\nu$.

By the previous theorem, the manifold $Q_\infty$ has a positive lower
bound to its injectivity radius since $\nu$ has bounded type.  If
$E_\infty$ represents the simply degenerate end of $Q_\infty$ for
which $\nu(E_\infty) = \nu$, then $E$ and $E_\infty$ represent ends of
two different manifolds with injectivity radius bounded below and the
same ending lamination.

Applying the main theorem of \cite{Minsky:ends}, or its generalization
\cite{Mosher:elc} if $N$ does not have a global lower bound to its
injectivity radius, the ends $E$ and $E_\infty$ are bi-Lipschitz
diffeomorphic.
\qed

\section{Grafting in degenerate ends}
\label{section:grafting}
In this section we describe a central construction of the paper.  The
{\em grafting} of a simply degenerate end, introduced as a technique in 
\cite{Bromberg:bers}, serves as the key to approximating degenerate
ends of complete hyperbolic 3-manifolds by cone-manifolds.

\bold{The grafting construction.}  
By Bonahon's theorem \cite{Bonahon:tame}, each manifold $M \in AH(S)$ is
homeomorphic to $S\times \reals$.  Given a particular choice of
homeomorphism $F \colon S\times \reals \to M$ each simple closed curve
$\alpha$ on $S$ has an associated embedded {\em positive grafting annulus} 
$$A^+ = F(\alpha \times [0 , \infty))$$ in $M$ and a {\em negative
grafting annulus} $A^- = F(\alpha \times (-\infty,0])$.

Consider the solid-torus cover $M_\zed =
\half^3/\langle  \alpha^*\rangle$  obtained as the quotient of
$\half^3$ by a representative of the conjugacy class of 
$F(\alpha \times \{0\})$ in
$\pi_1(M)$.  Then by the lifting theorem, the grafting annulus $A^+$
lifts to an annulus $A^+_\zed$ in the cover $M_\zed$.  

Let ${\rm Gr}^+(M,\alpha)$ denote the singular 3-manifold
obtained by isometrically  gluing the metric completions of 
$$M \setminus A^+ \ \ \ \ \text{and} \ \ \ \ M_\zed \setminus A^+_\zed$$
in the following way:
\begin{itemize}

\item[I.] For reference, choose an orientation on the curve $\alpha$.
Together with the product structure $F$, this orientation gives a
local ``left'' and ``right'' side in $M$ to the annulus $A^+$
corresponding to the left and right side of the curve $F(\alpha \times
\{t\})$ in $F(S \times
\{t\})$.  

\item[II.]
The metric completion of $M \setminus A^+$ contains two isometric
copies $\calA^l$ and $\calA^r$ of the annulus $A^+$ in its metric boundary
corresponding to the local left and right side of the annulus with
respect to the choice of orientation of $\alpha$.  Likewise, the
metric boundary of the metric completion of $M_\zed \setminus A^+_\zed$
contains the two isometric copies $\calA^l_\zed$ and $\calA^r_\zed$ of
$A^+_\zed$ corresponding to the local left and right side of
$A^+_\zed$ in $M_\zed$.

\makefig{Lifting the grafting
annulus.}{figure:lift}{\psfig{file=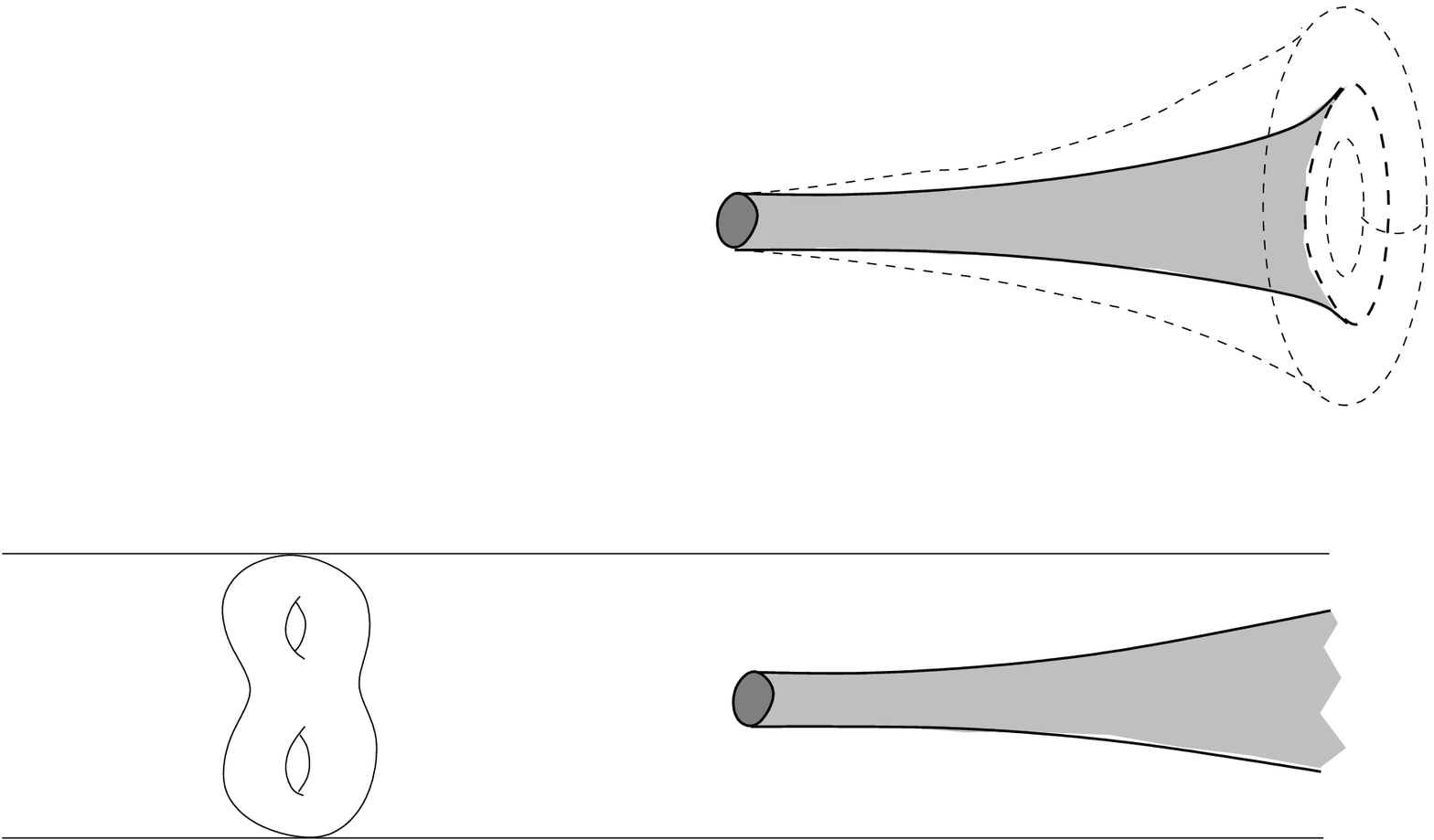,height=2.7in}}  

\item[III.]  The parameterization of $A^+$ by $F\vert_{\alpha \times [0,\infty)}$ induces
parameterizations 
$$F^l \colon \alpha \times [0,\infty) \to \calA^l 
\ \ \ \text{and} \ \ \ 
F^r \colon \alpha \times [0,\infty) \to \calA^r$$ 
of the annuli $\calA^l$ and $\calA^r$ and 
$$F_\zed^l \colon \alpha \times [0,\infty) \to \calA_\zed^l 
\ \ \ \text{and} \ \ \ 
F_\zed^r \colon \alpha \times [0,\infty) \to \calA_\zed^r$$
of the annuli
$\calA^l_\zed$ and $\calA^r_\zed$.
We obtain the
grafting ${\rm Gr}^+(M,\alpha)$ by identifying the metric completions
$$\closure{M \setminus A^+} \ \ \ \text{and} \ \ \ \closure{M_\zed
\setminus A^+_\zed}$$
by the mapping $\phi$ from the metric boundary of $M \setminus A^+$ to
the metric boundary of $M_\zed \setminus A_\zed^+$ determined by setting
$$\phi(F^l(x,t)) = F_\zed^r(x,t) \ \ \ \text{and} \ \ \ 
\phi(F^r(x,t)) = F_\zed^l(x,t).$$
(In a
neighborhood of $\alpha$, ${\rm Gr}^+(M,\alpha)$ is the two-fold
branched cover of $M$ branched along $\alpha$).

\item[IV.]  Since the gluing is isometric, the resulting grafted end
has a hyperbolic metric away from a singularity along the curve
$\alpha$.  When the curve $\alpha$ is a geodesic, the singularity
becomes a cone-type singularity, and ${\rm Gr}^+(M,\alpha)$ is a
hyperbolic cone-manifold homeomorphic to $S \times \reals$ with
cone-angle $4\pi$ at $\alpha$.

\makefig{The grafted end.}
{figure:graft}
{\psfig{file=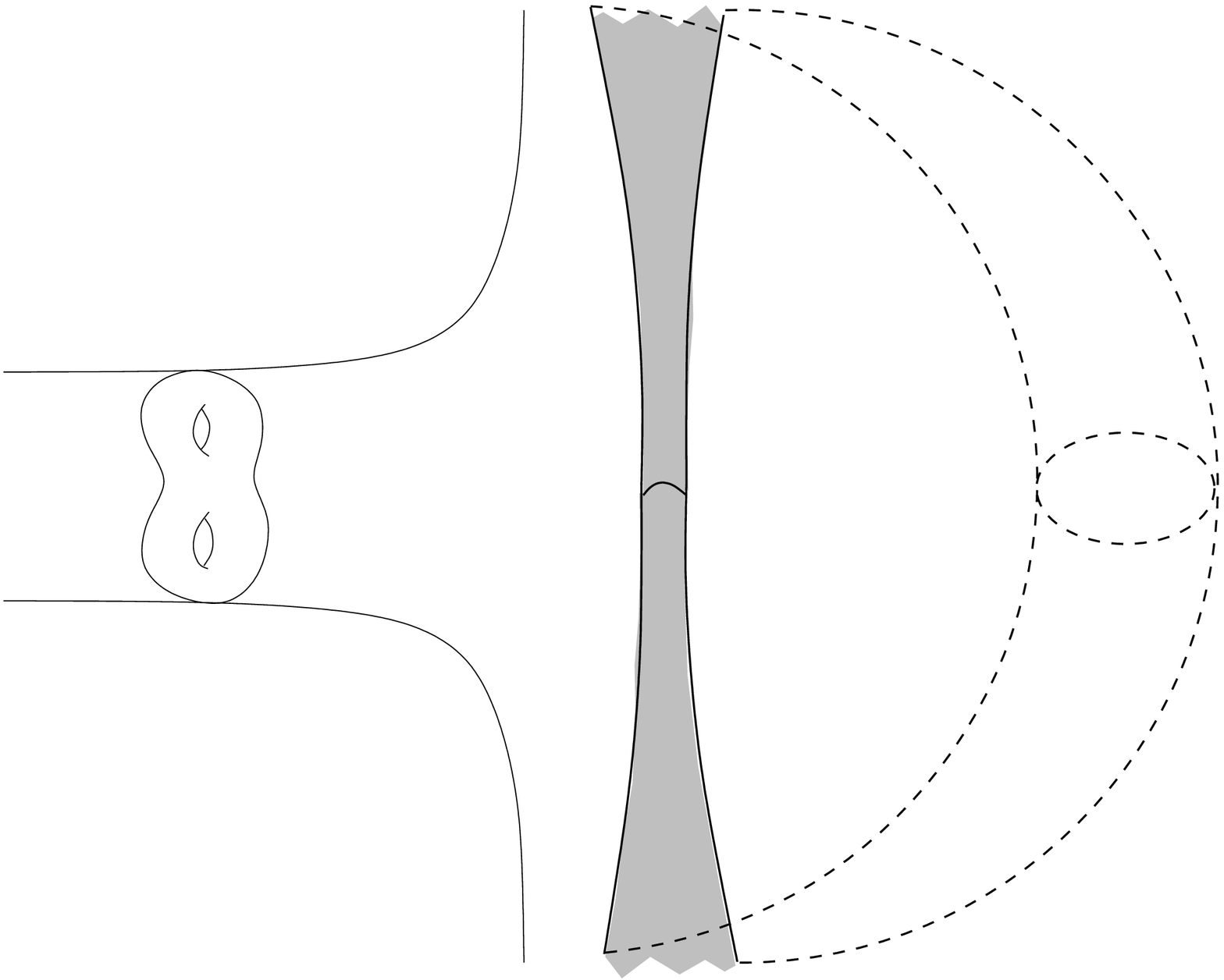,height=2.5in}}
\end{itemize}

By Otal's theorem (Theorem~\ref{theorem:otal}) any sufficiently short
geodesic in $M \in AH(S)$ is unknotted, guaranteeing that the grafting
construction may be applied to the geodesic itself.  We now prove that
grafting along a short geodesic always produces a geometrically finite end.
\begin{theorem}
If $(f \colon S \to M) \in AH(S)$ has no cusps and $\alpha$ is an
essential simple closed curve in $S$ for which $\ell_M (\alpha^*) <
\ell_{\rm knot}$ then the positive end of the hyperbolic
cone-manifold ${\rm Gr}^+(M,\alpha^*)$ is geometrically finite.
\label{theorem:cone:density}
\end{theorem}

\bold{Proof:} 
Applying Otal's theorem (Theorem~\ref{theorem:otal}), $\alpha$ is
isotopic into a level surface for 
any product structure on $M$, so we choose a homeomorphism
$F \colon S \times \reals \to M$ so that 
\begin{enumerate}
\item $F(S \times \{0\})$ is homotopic to
$f$, and
\item $F\vert_{S \times \{0\}}$ {\em realizes} $\alpha$: i.e.,
$F(\alpha \times \{0\}) = \alpha^*$.
\end{enumerate}

Let $E^+$ be the positive end of $M$.  Let $M^c = {\rm
Gr}^+(M,\alpha^*)$ and let $E^c$ denote the positive end of $M^c$.
Arguing by contradiction, assume the end $E^c$ is not geometrically
finite.  Section~\ref{section:geom:finite} guarantees, then, that
there are simple closed curves $\gamma_k$ on $S$ whose geodesic
representatives in $E^c$ eventually lie outside of every compact
subset of $E^c$.

We choose a particular exhaustion of $E^c$ by compact submanifolds
that is adapted to an exhaustion of the original end $E^+$ as
follows:
\begin{enumerate}
\item Let $K_j$ be an exhaustion of $E^+$ by the compact submanifolds
$$K_j = F(S \times [0,j]).$$
\item Let $\hat{K_j}$ be the lift of $K_j \setminus A$ to $E^c$ for which the
restriction of $\pi$ to $\hat{K_j}$ is an isometric embedding.
\item Extend $\hat{K_j}$ to a compact subset $K_j^c$ by taking 
the union of $\hat{K_j}$ with an exhaustion of $M_\zed$ by solid tori
as follows: let $F_\zed$ be the lift of $F_{\alpha \times [0,\infty)}$
to $M_\zed$, and let $V_j$ be an exhaustion of $M_\zed$ by closed
solid tori so that $V_j$ intersects 
$A^+_\zed$ in $F_\zed(\alpha \times [0,j])$.  
Then $$K_j^c = \hat{K_j} \cup V_j$$ exhausts the
end $E^c$.
\end{enumerate}

Let $\gamma_j \subset \{ \gamma_k\}$ be a subsequence for which
$$\gamma_j^* \subset E^c \setminus K_j^c.$$

We claim there is a compact subset $K$ of $M$ so that the projections
$\pi(\gamma_j^*)$ of $\gamma_j^*$ to $M$ all intersect $K$.  Consider
the (unique) component $S_\zed$ of the lifts $S$ to the solid torus
$M_\zed$ for which $\pi_1(S_\zed) = \zed$, i.e. $S_\zed$ is the
annular lift of $S$ to $M_\zed$ that contains the curve $\alpha$.  The
properly embedded annulus $S_\zed$ separates $M_\zed$ into two pieces,
one covering the component of $M \setminus S$ containing $S_0$, and
the other covering the noncompact portion of $E^+ \setminus S$. 

Since $\alpha \subset K_j^c$ for all $j$ sufficiently large, we
may throw away a finite number of $\gamma_j$ to guarantee that 
$\alpha \not= \gamma_j$ for all $j$.
Thus,  the geodesic $\gamma_j^*$ intersects $S_\zed$ if and
only if we have
$$i(\gamma_j, \alpha) \not= 0.$$

The geodesic $\gamma_j^*$ projects isometrically by
$\pi$ to the geodesic representative of $\gamma_j$ in $M$; for
convenience, we denote the latter by $\pi(\gamma_j^*)$.  
Let $X_j$ be a pleated surface realizing $\gamma_j$ in $M$ with the
property that if $i(\gamma_j, \alpha) = 0$ then $X_j$ realizes
$\alpha$ as well.

If $\gamma_j^*$ intersects $S_\zed$ in $M^c$, then the
geodesic $\pi(\gamma_j^*)$ intersects $S$ in $M$, so the pleated
surface $X_j$ intersects $S$.  If, on the other hand, $\gamma_j^*$
does not intersect $S_\zed$, then $X_j$ realizes $\alpha$,
so $X_j$ also intersects $S$.  By Theorem~\ref{theorem:diameter},
there is a compact subset 
$K \subset M$ so that we 
have $$X_j \subset K$$ for all $j$.  In particular, it follows that we
have $$\pi(\gamma_j^*) \subset K$$ for all $j$.

Note, however, that there is a $j'>j$ so that $\gamma_j^*$
intersects $\hat{K_{j'}} \setminus \hat{K_j}$ since otherwise $\gamma_j^*$ would 
lie entirely in $M_\zed \setminus A^+_\zed$ which would imply
$\gamma_j$ is isotopic to $\alpha$.  Choosing $j$ sufficiently large 
to guarantee that 
$$K \subset K_j,$$ then, we obtain a contradiction, since
$$\pi( \hat{K_{j'}} \setminus \hat{K_j} ) \cap K = \nullset.$$
We conclude that the end $E^c$ is geometrically finite.
\qed

We introduce one further piece of notation for later use.  If $\alpha$
and $\beta$ represent simple closed curves whose geodesic
representatives lie in $E^-$ and $E^+$ respectively, we can perform
grafting of $E^-$ along the negative grafting annulus for $\alpha$ in
$M$ and grafting of $E^+$ along the positive grafting annulus for $\beta$ in
$M$ simultaneously.  We denote by ${\rm Gr}^\pm(M,\alpha,\beta)$ this
{\em simultaneous grafting along $\alpha$ and $\beta$}.


\section{Geometric inflexibility of cone-deformations}
\label{section:cone}
\newcommand{\sym}{{\rm sym}}
\newcommand{\curl}{\operatorname{curl}}
\newcommand{\cL}{{\calL}}
\newcommand{\dev}{{\sf Dev}}
\newcommand{\ray}{{\sf r}}

In this section we establish the estimates on cone-deformations
necessary to obtain geometric control away from the cone-singularity.
The main result of this section, Theorem~\ref{theorem:bilip}, harnesses the
deformation theory of cone-manifolds originally developed in
\cite{Hodgson:Kerckhoff:rigidity} and \cite{Bromberg:thesis},
and improved upon in
\cite{Hodgson:Kerckhoff:bounds,Hodgson:Kerckhoff:shape} and 
\cite{Bromberg:Schwarzian}  to obtain the necessary
control.  

Theorem~\ref{theorem:bilip} is a geometric formulation of
similar analytic results in \cite{Bromberg:Schwarzian} that control
the change of projective structure associated to an end of
cone-manifold $M$ under a deformation that changes the cone-angle.  Here,
we have 
replaced control over the change in 
projective structure, which suffices for applications in 
\cite{Bromberg:bers}, with bi-Lipschitz control over the metric
on $M$ itself.  We refer the reader to
\cite{Hodgson:Kerckhoff:harmonic} for an expository account of the
necessary recent developments in the cone deformation theory.

\smallskip

Let $E$ be an end of a geometrically finite cone-manifold that is cut
off by a surface $S$ with genus $\ge 2$. Then the hyperbolic structure
on $E = S \times \reals^+$ naturally extends to a conformal structure
on $S \times \{\infty\}$. The hyperbolic structure on $E$ is locally
modeled on $\half^3$, while the conformal structure is modelled on
$\chat$. 

More concretely, $\half^3$ is compactified by $\chat$, and $\PSL_2\cx$
acts continuously on the compactification as hyperbolic isometries of
$\half^3$ and as projective transformations of $\chat$. Then $S \times
 (0, \infty]$ has an atlas of charts to $\half^3 \cup \chat$ with
transition maps restrictions of elements of $\PSL_2\cx$. These charts
will map points in $S \times \reals^+$ to $\half^3$ and points in $S
\times \{\infty\}$ to $\chat$. Restricted to $S \times \reals^+$ the
charts will be an atlas for the end $E$, and on $S \times \{\infty\}$
the charts will define a conformal structure on $S$. We refer to $S$
with this conformal structure as the {\em conformal boundary} of $E$.

The following theorem appears in \cite{Bromberg:Schwarzian}.
\begin{theorem}
Given $\alpha >0$ there exists $\ell >0$ such that the following
holds.  Let $M_\alpha$ be a geometrically-finite hyperbolic
cone-manifold with no rank-1 cusps, cone-singularity $\calC$ and cone-angle
$\alpha$ at each component $c \subset \calC$.  If the tube-radius $R$ about
each component $c$ in $\calC$ is at least 
$\sinh^{-1}(\sqrt{2})$  and the total length
$\ell_{M_\alpha}(\calC)$ of $\calC$ in $M_\alpha$ satisfies
$\ell_{M_\alpha}(\calC) < \ell$ then there is a one parameter family
$M_t$ of cone-manifolds with fixed conformal boundary and cone-angle $t \in [0,\alpha]$ at each $c
\subset \calC$.
\label{theorem:one:parameter}
\end{theorem}

The main result of this section allows us to control the geometric
effect of the one-parameter cone-deformation when the cone-singularity
is sufficiently short.  
\begin{theorem}{\sc The Drilling Theorem}
Given $\alpha>0$, $L >1$, there exists $\ell>0$ so that the
following holds.  If $M_\alpha$ is a hyperbolic cone-manifold
satisfying the hypotheses of the previous theorem, and $M_t$ the
corresponding one-parameter family of cone-manifolds with $t \in [0,
\alpha]$, then if $\ell_{M_\alpha}(\calC) < \ell$ there is for
each $t$ a standard neighborhood $\tube_t(\calC)$ of the cone locus
$\calC$ and an
$L$-bi-Lipschitz diffeomorphisms of pairs
$$h_t \colon 
(M_\alpha \setminus \tube_\alpha(\calC), \bdry \tube_\alpha(\calC)) 
\to 
(M_t \setminus \tube_t(\calC), \bdry \tube_t(\calC))$$
so that $h_t$ extends to a homeomorphism
$\overline{h_t} \colon M_\alpha \to M_t$ for each $t \in (0,\alpha]$.
\label{theorem:bilip}
\end{theorem}

As we will see, the standard neighborhood $\tube_t(\calC)$ will be a
components $\tube_t^{\epsilon}(\calC)$ of the {\em Margulis
$\epsilon$-thin part} of $M_t$ containing 
$\calC$.  In fact, for each $\epsilon >0$ less 
than the appropriate {\em Margulis constant} there is an $\ell$ 
and $h_t$ satisfying the theorem for $\tube_t(\calC) =
\tube_t^\epsilon(\calC)$.

\bold{Background.}
Before proving Theorem~\ref{theorem:bilip} we review necessary
background.
Let $N$ be a 3-manifold with boundary (we allow $N$ to be
non-compact).
Let $g$ be a hyperbolic metric on the interior of $N$ that extends to 
a conformal structure on each component of $\bdry N$; here, the
metric $g$ need not be complete, but the conformal structures
compactify the ends where the metric is complete.

Let $\wt N$ denote the universal cover of $N$ and let $\pi \colon \wt
N \to N$ denote the covering projection.
Then $g$ lifts to a metric $\wt g$ on the universal cover
$\interior(\wt N)$ of $\interior(N)$, and the
conformal structures on $\bdry N$ lift to conformal structures on 
$\bdry \wt N$.  
There is a map $$\dev \colon \wt N \to \half^3 \cup \chat$$
that is a local isometry on $\interior(N)$ and locally
conformal on $\bdry \wt N$.  

Furthermore, there is a representation
$$\rho \colon \pi_1(N) \to \PSL_2(\cx)$$ with the property that
\begin{equation}
\dev(\gamma(p)) = \rho(\gamma) \dev(p)
\label{equation:develop}
\end{equation}
for each $p \in \wt N$ and each deck-transformation $\gamma \in
\pi_1(N)$.  

The map $\dev$ is called the {\em developing map} for the metric
$g$ and is determined up to post-composition with elements of
$\PSL_2(\cx)$ acting on $\half^3 \cup \chat$.  Changing the developing
map by postcomposition changes the corresponding representation by
conjugation in $\PSL_2(\cx)$.

A smooth family $g_t$ of such metrics on $N$ determines a smooth
family of developing maps
$$\dev_t \colon \wt N \to \half^3 \cup \chat.$$
The developing maps $\dev_t$ determine a time-dependent vector field
$v_t$ on $\wt N$, where
$v_t (p)$ is the pull-back by $\dev_t$ of the tangent vector to the
path $\dev_t(p)$ at time $t$; i.e.
$$v_t(p) = \dev^*_t \left( \frac{d \dev_t}{dt}(p) \right).$$
We call $v_t$ the {\em derivative} of the family of developing maps $\dev_t$.

A {\em Killing field} on a Riemannian manifold $N$ is a vector field
whose local flow is an isometry for the Reimannian metric on $N$.
 By differentiating \eqref{equation:develop} we have that for
each $t$ and for any
$\gamma \in \pi_1(N)$, the vector field
$$v_t - \gamma^* v_t$$
is a Killing field on $\interior(\wt N)$ for the Riemannian metric $g_t$.
A vector field with this automorphic property 
is called {\em automorphic for the metric $g_t$},
or {\em $g_t$-automorphic}.  Note that a $g_t$-automorphic vector
field $v_t$ need not arise as the derivative of a family of developing
maps. 

Let $\hat g_t$ be a family of metrics on $N$ for which there are
diffeomorphisms $f_t \colon N \to N$ isotopic to the identity
satisfying $(f_t)^* g_t = \hat g_t$, i.e.  $$f_t \colon (N,\hat g_t)
\to (N,g_t)$$ is an isometry.  Then there will be a corresponding
family of developing maps $\hat \dev_t$ with derivative $\hat v_t$ a
time-dependent vector field on $\wt N$.

The lifts $\wt f_t \colon \wt N \to \wt N $ of $f_t$ to the universal
cover allow us to compare $v_t$ and $\hat v_t$: since the maps $f_t$
are isometries from $(N,\hat g_t)$ to $(N,g_t)$ one sees that the
difference $$(\wt f_t)_* \hat v_t - v_t$$ restricts to the sum of a
$\pi_1(N)$-{\em equivariant} vector field and a Killing field on
$\interior(\wt N)$.
In fact,
$\hat \dev_t$ can be chosen so that $(\wt f_t)_* \hat v_t - v_t$ is a
$\pi_1(N)$-equivariant vector field: one simply needs to alter $\hat
\dev_t$ by postcomposition with a family of elements in $\PSL_2(\cx)$.

The derivative of the family of developing maps $\dev_t$ is a
$g_t$-automorphic vector field; we now seek to integrate a
$g_t$-automorphic vector field $v_t$ to obtain a family of developing
maps $\dev_t$, reversing the above process.

\begin{theorem}
Let $w_t$ be a smoothly varying uniquely integrable family of
$g_t$-automorphic vector fields on $\wt N$ tangent to the boundary
$\bdry \wt N$ such that $w_t - v_t$ is equivariant.  For any subset $U
\subset N$ contained in a compact subset of $N$, there exists a family
of metrics $\hat g_t$ on $N$, developing maps $\hat \dev_t$ for $\hat
g_t$, $\hat g_t$-automorphic vector fields $\hat v_t$ on $\wt N$, and
diffeomorphisms $f_t \colon N \to N$ isotopic to the identity 
so that the following holds:
\begin{enumerate}
\item $(f_t)^* g_t = \hat g_t$,
\item $\hat v_t$ is the derivative of the developing maps $\hat \dev_t$,
and
\item $(f_t)_* \hat v_t = w_t$ on 
$\pi^{-1}(U)$.
\end{enumerate}
\label{theorem:develop}
\end{theorem}

\bold{Proof:}  We first prove that each $p \in N$ 
has a neighborhood $U$ such that the
theorem holds on $U$ for $t$ near $0$.   

Let $\wt p \in \pi^{-1}(p)$ and choose nested neighborhoods $V
\subset V' \subset V''$ of $\wt p$
such that both $\dev_0$ and $\pi$ restricted to $V''$ are embeddings.
Then there exists an $\epsilon' >0$ such that
for $|t| < \epsilon'$ the image $\dev_t(V'')$ contains
$\dev_0(V')$. We can then choose an $\epsilon''$ with $0< \epsilon''
\le \epsilon'$ such that  
for $|t| < \epsilon''$ there
exists a flow
$$\phi_t \colon \dev_0(V) \to \half^3 \cup \chat$$
so that $\phi_t$ is the flow of the time-dependent vector fields
$(\dev_t)_* w_t$.

Given $q \in V$ define 
$$\hat \dev_t(q) = \phi_t \circ \dev_0(q).$$
Then we claim $\hat \dev_t$ can be extended to a developing map on all of
$\wt N$ for each $t$.

To see this, we first note that for $|t| < \epsilon''$, we have
$$\hat \dev_t (V) \subset \dev_t(V'')$$
so we may define 
 an embedding
$f_t \colon V \to \wt N$ by 
$$h_t = \dev_t^{-1} \circ \hat \dev_t.$$
Since $V''$ is disjoint from its translates, there exists a
family of embeddings 
$f_t \colon \pi(V) \to N$ 
so that $$\wt f_t \vert_V = h_t$$
on $V$.  
We extend $f_t$ to a smooth family of diffeomorphism of all of $N$ so that 
$f_0$ is the identity.  
Then by setting $U = \pi(V)$ and 
$$
\hat \dev_t = \dev_t \compos \wt f_t 
$$ 
we obtain the desired extension.

\smallskip

Now let $U \subset N$ be any subset of $N$ contained in a compact
subset $K$ of $N$.  We again establish the theorem for $U$ and for $t$
near $0$. By the above, and the compactness of $K$, there exists a
finite collection of open sets $U_i \subset N$ and an $\epsilon >0$ so
that $\cup_i U_i$ covers $K$ and the theorem holds on each $U_i$ for
$|t|< \epsilon$.  Let $\dev^i$ be the resulting developing maps for each
$U_i$.

By the uniqueness of flows, we may define $\hat \dev_t$ by 
$$\hat \dev_t \vert_{\pi^{-1} (U_i)} = \dev^i_t\vert_{\pi^{-1}(U_i)},$$ 
and we have the theorem for $U$ and all $|t| < \epsilon$.

It follows that we may now define $\hat \dev_t$ on some open interval
$(a,b)$.  Applying the theorem at $t = b$ we have  corresponding
triples
$(g_t',\dev_t',v_t')$ satisfying the conclusions of the theorem on $U$
for $t \in (b- \epsilon', b+ \epsilon')$ for some $\epsilon' >0$.
Let $f_t' \colon (N,g_t') \to (N,g_t)$ be the corresponding
diffeomorphisms.

The developing maps $\dev_t'$ satisfy
$$\dev_t' = \dev_t \compos \wt f_t'.$$
Thus, we have
$\dev_t = \dev_t' \compos (\wt f_t')^{-1}$ and therefore setting
$$\hat \dev_t = \dev_t' \compos (\wt f_t')^{-1} \compos \wt f_t$$
extends $\hat \dev_t$ over a neighborhood of $t=b$.  
Arguing similarly for $t=a$, the set of $t$ values on which $\hat \dev_t$
may be defined is open, closed and non-empty, and therefore $\hat \dev_t$
can be defined for all $t$.  
The proof is complete.
\qed

Let $g_t$ be a smooth family of Riemannian metrics on $N$. We define
vector valued 1-forms $\eta_t$ by the formula 
$$\frac{dg_t(x,y)}{dt} = 2g_t(x, \eta_t(y)).$$
The symmetry of $g_t$ implies that $\eta_t$ is self-adjoint. We define
a pointwise norm of $\eta_t$ by choosing an orthonormal frame field
$\{e_1, e_2, e_3\}$ for the $g_t$-metric and setting 
$$\|\eta_t\|^2 = \sum_{i,j} g_t(\eta_t(e_i), \eta_t(e_j)).$$
Note that
$$g_t(x, \eta_t(x)) \leq \|\eta_t\| g_t(x,x).$$

Given two metrics $g$ and $\hat{g}$ we define the biLipschitz constant
at each point $p \in N$ by 
$${\rm bilip}_p(g, \hat{g}) = \sup \left\{
K \geq 1\Big\vert \frac{1}{K} \le \sqrt{\frac{\hat{g}(x,x)}{g(x,x)}}
\le K \mbox{\ for all $x \in T_p N$, $x \neq 0$} \right\}.$$ A
bound on $\|\eta_t\|$ for all $t\in [0,a]$ gives a bound on ${\rm
bilip}_p(g_0, g_a)$. In particular $$\left| \frac{dg_t(x,x)}{dt}
\right| \leq 2\|\eta_t\| g_t(x,x)$$ and integrating we have $$g_a(x,x)
\leq e^{2Ka} g_0(x,x)$$ if $\|\eta_t\| \leq K$ for all $t \in
 [0,a]$. This implies that $${\rm bilip}_p(g_0, g_a) \leq e^{aK}.$$

The families of metrics we will examine will always be the pullback of
some fixed metric $g$ by the flow $\phi_t$ of a time dependent vector
field $v_t$. In this 
case we can relate $\eta_t$ to the covariant derivative of $v_t$. More
precisely if $g_t = \phi_t^* g$ then $$\eta_t = \sym \nabla^t v_t$$
where $\nabla^t$ is the Riemannian connection for the $g_t$-metric and
$\sym \nabla^t$ is the symetric part of the covariant derivative. This
follows from the fact that 
$$\frac{dg_t(x,y)}{dt} = {\cal L}_{v_t} g_t(x,y) = g(\nabla^t_x v_t,
y) + g(x, \nabla^t_y v_t) = 2 g(x, \sym \nabla^t_y v_t).$$  

Our vector fields will also be divergence free and harmonic. For our
purposes $v$ is harmonic if it is divergence free and $\curl \curl v =
-v$. Note that the our $\curl$ is half the usual $\curl$ and is chosen
to agree with the definition given in
\cite{Hodgson:Kerckhoff:rigidity}. We also refer there for motivation for
this definition of harmonic. Note that $\curl v$ will also be a
divergence free, harmonic vector field. 

We use $\nabla^t$ to define an operator $D_t$ on the space of vector
valued $k$-forms by the formula $$D_t = \sum_i \omega^i \wedge \nabla^t_{e_i}$$
where the $e_i$ are an orthonormal frame field with coframe
$\omega^i$. The formal adjoint of $D_t$ is then $$D_t^* = \sum_i
i(e_i) \nabla_{e_i}$$ where $i(e_i)$ is contraction.

Let $w_t = \curl v_t$.
In section 2 of \cite{Hodgson:Kerckhoff:rigidity} it is shown that
$$\sym \nabla^t w_t = *D_t \eta_t = \beta_t$$ and $D^*_t \eta_t = 0$. 
Bounds on the
norms of $\eta_t$ and $\beta_t$ will allow us to control the geodesic
cuvature of a smooth curve in $N$.

\begin{prop}
\label{curvature}
Let $\gamma(s)$ be smooth curve in $N$ and let $C(t)$ be the geodesic
curvature of $\gamma$ at $\gamma(0) = p$ in the $g_t$ metric. For each
$\epsilon>0$ there exists a $K>0$ depending only on $\epsilon$, $a$
and $C(0)$ such that $|C(a) - C(0)| \le \epsilon$ if $\|\eta_t(p)\|
\leq K$, $\|\beta_t(p)\| \le K$ and $D_t^* \eta_t = 0$ for all $t \in
[0,a]$.
\end{prop}

\bold{Proof:}  We assume that $\gamma(s)$ is a unit speed
parameterization in the $g_0$-metric. In the $g_t$-metric we
reparameterize such that $\gamma_t(s) = \gamma(h_t(s))$ is a unit
speed parameterization. Let 
$$V(t) = \nabla^t_{\gamma'_t} \gamma'_t(0).$$
Then
$$C(t)^2 = g_t(V(t), V(t)).$$
Differentiating we have
\begin{equation}
\label{cderform}
C(t) C'(t) = g_t(V(t), V'(t)) + g_t(V(t), \eta_t(V(t))).
\end{equation}
The result will follow if we can bound
$C'(t)$; we accomplish this via a calculation in local coordinates. We
choose our coordinates so as to bound the derivative at $t=0$.

To write the various tensors in local
coordinates $(x_1, x_2, x_3)$ we let $e_i =
\frac{\del}{\del x_i}$, define functions $g_{ij}$ by $g_t(e_i,
e_j) = g_{ij}(t)$, and let the $g^{ij}$ be chosen such that
$(g_{ij})(g^{ij}) = \id$. We similarly define $\eta_i^j$ by the
formula $\eta_t(e_i) = \sum_j \eta_i^j(t) e_j$. The $\beta_i^j$ are
defined in the same way. Recall that the Christoffel symbols
$\Gamma_{ij}^k$ satisfy the formula $$\Gamma_{ij}^k = \frac{1}{2}
\sum_{m} \left(\frac{\del g_{im}}{\del x_j} + \frac{\del g_{jm}}{\del
x_i} - \frac{\del g_{ij}}{\del x_m}
\right) g^{mk}.$$ 
We can choose our local coordinates such that at a point $p$ in $N$ we
have $g_{ij}(0) = \delta_i^j$ and $\Gamma_{ij}^k (0) = 0$. Note that
$\frac{1}{2} \frac{dg_{ij}}{dt} = \sum_k \eta_i^k g_{kj}$ so 
with this choice of coordinates we have 
$\eta_i^j(0) = \frac{1}{2}\frac{d
g_{ij}}{dt}(0)$ and
\begin{equation}
\label{christ}
\frac{d \Gamma_{ij}^k}{dt}(0) = \frac{\del \eta_i^k}{\del x_j}(0) +
\frac{\del \eta_k^j}{\del x_i}(0) - \frac{\del \eta_i^j}{\del x_k}(0)
\end{equation}
at $p$.

We need to write the $\beta_i^j$ as derivatives of the $\eta_i^j$. To
do this we note that by definition 
$$\eta_t = \sum_{i,j} \eta_i^j(t) e_j 
\otimes \omega^i.$$
Then direct calculation gives
$$D_0 \eta_0 = \sum_{i,j,k} \frac{\del \eta_i^j}{\del x_k}(0) e_j
\otimes \omega^k \wedge \omega^i$$ 
at $p$. This implies that
\begin{equation}
\label{Deta}
\beta_i^j(0) = \frac{\del \eta_{i+2}^j}{\del x_{i+1}}(0) - \frac{\del
\eta_{i+1}^j}{\del x_{i+2}}(0) 
\end{equation}
at $p$, where on the right hand side of this formula the coefficients
are measured mod 3. From this we see that a bound on $\|\beta\|$ gives
a bound on the difference $$\frac{\del \eta^j_i}{\del x_k}(0) -
\frac{\del \eta^j_k}{\del x_i}(0).$$

Another direct calculation in local coordinates gives
\begin{equation}
\label{Dstareta}
D_0^* \eta_0 = \sum_{i,j} \frac{\del \eta_i^j}{\del x_i}(0) e_j.
\end{equation}
Therefore, if 
$D_0^* \eta_0 = 0$ we have $\frac{\del \eta_i^j}{\del x_i}(0) = 0$.

Combining (\ref{christ}), (\ref{Deta}) and (\ref{Dstareta}) we have
\begin{eqnarray}
\frac{d\Gamma_{11}^1}{dt}(0) & = & 0 \nonumber \\
\frac{d\Gamma_{11}^2}{dt}(0) & = & \beta_3^1(0) \label{christbeta} \\
\frac{d\Gamma_{11}^3}{dt}(0) & = & -\beta_2^1(0) \nonumber.
\end{eqnarray}

We can now return to our smooth curve
$\gamma(s)= (x_1(s), x_2(s), x_3(s))$. Assume that $\gamma(0) = p$ and
$\gamma'(0) = e_1$. Define $h_1(t) = \frac{dh_t}{ds}(0)$ and
$h_2(t) = \frac{d^2 h_t}{ds^2}(0)$.
By our choice of $\gamma$ we have
$$V(t) = h_1(t)^2\sum_i \left(\frac{d^2
x_i}{ds^2}(0) + \Gamma_{11}^i(t) \right) e_i + h_2(t)e_1.$$
To bound $C'(0)$ we need to bound $g_0(V(0), V'(0))$ and $g_0(V(0),
\eta_0(V(0)))$. For the second term we have 
$$|g_0(V(0), \eta_0(V(0)))| \leq C^2(0)\|\eta_0\|.$$
To bound the first term we note that $V(0)$ is perpendicular to
$\gamma'(0)$. Therefore, if we let $V_{\bot}(t)$ be the sum of the
$e_2$ and $e_3$ terms of $V(t)$, we have $g_0(V(0), V'(0)) = g_0(V(0),
V'_{\bot}(0))$. 
Differentiating we have
$$V'_\bot(0) =  \sum_{i=2,3} \left(2h'_1(0) \frac{d^2 x_i}{ds^2}(0) + \frac{d
\Gamma_{11}^i}{dt} (0) \right) e_i.$$
We need to bound $h'_1(0)$. This is accomplished by
differentiating the formula
$$g_t(\gamma'_t(s), \gamma'_t(s)) = 1.$$
Note that when $s=0$ the formula becomes
$$h_1(t)^2 g_{11}(t) = 1$$
and differentiating with respect to $t$ yields $h'_1(0) = - \eta_1^1(0)$.
To bound the derivative of the Christoffel symbols we use
\eqref{christbeta}. Therefore 
$$\left|V'_\bot(0) \right|_0 \le  2C(0)\|\eta_0\| + \|\beta_0\|$$
which in turn implies that
$$|g_0(V(0), V'(0))| \leq C(0)(2C(0)\|\eta_0\| + \|\beta_0\|)$$
and
$$|C'(0)| \leq 3C(0) \|\eta_0\| + \|\beta_0\|.$$
Since we could choose coordinates to calculate the derivative for any
$t$ we have 
\begin{equation}
\label{curvatureder}
|C'(t)| \leq 3C(t) \|\eta_t\| + \|\beta_t\|.
\end{equation}
Therefore if we choose $K$ small enough if $\|\eta_t(p)\|$ and
$\|\beta_t(p)\|$ are less than $K$, then integrating
(\ref{curvatureder}) implies that $|C(a) - C(0)| \le \epsilon$. \qed 

We remark that a similar statement holds for subsurfaces of $N$. In
particular if $\gamma$ lies on a subsurface $S$ then the metrics $g_t$
induce a metric on $S$. If we use this induced metric on $S$ to
measure the geodesic curvature of $\gamma$ then the conclusion of
Proposition \ref{curvature} still holds. 

For a vector field $v$ the symetric traceless part of $\nabla v$ is
the {\em strain} of $v$ and measures the conformal distortion of the
metric pulled back by the flow of $v$. If $v$ is divergence free then
$\nabla v$ will be traceless so $\eta = \sym \nabla v$ is a strain
field. If $v$ is harmonic we say that $\eta$ is also harmonic. To
prove Theorem~\ref{theorem:bilip} we need the following mean-value 
inequality for harmonic strain fields (the result is due to Hodgson
and Kerkchoff,  see \cite[Thm. 9.9]{Bromberg:Schwarzian} for an
exposition). 

\begin{theorem}  Let $\eta$ be a harmonic strain field on a ball $B_R$
of radius $R$ centered at $p$.  Then we have
$$\| \eta(p) \| \le
\frac{ 3 \sqrt{2 \vol(B_R)}}{4 \pi f(R) } 
\sqrt{
\int_{B_R} \| \eta \|^2 dV}
$$
where 
$$f(R) = \cosh(R) \sin (\sqrt{2} R) - \sqrt{2} \sinh (R) \cos
(\sqrt{2} R)$$ and $R < \pi/\sqrt{2}$.
\label{theorem:mvt}
\end{theorem}

We now return to the concrete situation of interest: we assume given
triples $(g_t,\dev_t,v_t)$, a family of cone-metrics, developing maps,
and $g_t$-automorphic vector fields $v_t$ so that $v_t$ is the
derivative of $\dev_t$, where $t \in [0,\alpha]$ denotes the cone-angle
at the cone-singularities of $g_t$.  To apply
Theorem~\ref{theorem:mvt}, we will invoke the following Hodge theorem
of Hodgson and Kerckhoff \cite{Hodgson:Kerckhoff:rigidity} and its
generalization \cite{Bromberg:thesis} to the geometrically finite
setting.

\begin{theorem}{\sc Hodge Theorem}
Given the triple $(g_t, \dev_t, v_t)$ 
there exists a smooth, time-dependent, divergence-free, harmonic,
$g_t$-automorphic vector field $w_t$
so that for each $t\in [0,\alpha]$ we have 
\begin{enumerate}
\item $w_t$ is tangent to $\bdry \wt N$,
\item the restriction of $w_t$ to $\bdry \wt N$ is conformal, and 
\item $w_t - v_t$ is an equivariant vector field.
\end{enumerate}
\end{theorem}

Combined with Theorem~\ref{theorem:develop} the Hodge Theorem has the
following corollary. 

\begin{cor}
\label{nicemetric}
There exists a one-parameter family of cone-metrics $g_t$ on $N$ such
that $M_t = (N, g_t)$ and $\eta_t$ is a harmonic strain field outside
a small tubular neighborhood of the cone-singularity and the rank-2
cusps. 
\end{cor}

Below, we will estimate the $L^2$-norm of $\eta_t = \sym \nabla w_t$
outside of a tubular neighborhood of a short cone-singularity.  
This, together with Theorem~\ref{theorem:mvt} will give us the
necessary control metrics $g_t$ outside of the thin part. 
Before obtaining this control, we must normalize the picture in a
neighborhood of the cone-singularity.

In general, the Margulis lemma does not apply to cone-manifolds.  If,
however, there is a uniform lower bound $R$ to the tube-radius of each
component of the cone-singularity and an upper bound $\alpha$ on all
cone-angles, a thick-thin decomposition exists exactly analogous to
that of the smooth hyperbolic setting (see
\cite{Hodgson:Kerckhoff:bounds,Bromberg:Schwarzian}).  In patricular
there exists and $\epsilon_{R, \alpha}$ such that the
$\epsilon_{R,\alpha}$-thin part $M^{\le \epsilon_{R,\alpha}}$ of a
hyperbolic cone-manifold $M$ consists of tubes about short geodesics
(including the cone singularity) and cusps.   

In our situation, we have assumed that the cone singularity of
$M_\alpha$ has tube-radius at least $\sinh^{-1} \sqrt{2}$. It is shown
in
\cite{Hodgson:Kerckhoff:bounds}
that this tube radius will not decrease as the cone-angle
decreases. Therefore we fix 
$$\varepsilon = \epsilon_{\sinh^{-1}
\sqrt{2}, \alpha}.$$

Given a non-parabolic homotopy class $[\gamma]$ of a closed curve
$\gamma$ in $M_\alpha$, 
it will be convenient to consider the family of embedded tubes with
core the geodesic representative of $\gamma$ as the cone-angle varies.
For this, we use the following notation: if $\gamma$ is homotopically
non-trivial closed curve in $N$ 
with $\ell_{M_t}(\gamma) \le \epsilon \le \varepsilon$, we denote by 
$\tube^{\epsilon}_t(\gamma)$ is the component of $M^{< \epsilon}_t$
that contains the geodesic representative of $\gamma$ in the
$g_t$-metric. We will often need to refer to the union of the tubes 
about the cone-singularity $\calC$. For this reason we set 
$$\tube^\epsilon_t(\calC) = \underset{c \in \calC}{\cup}
\tube_t^\epsilon(c).$$
Occasionally we will make statements about a generic Margulis tube
without reference to a particular tube in the cone manifolds $M_t$. We
simply refer to such a generic $\epsilon$-Margulis tube as
$\tube^\epsilon$.   

\begin{theorem}  Given $\epsilon >0$, there is an $\ell >0$ and $K>0$
such that if $\ell_{M_t}(\calC) \le \ell $  then we have the $L^2$-bound  
$$\int_{M_t \setminus \tube^\varepsilon_t(\calC)} \|  \eta_t \|^2  +
\|*D_t \eta_t\|^2 \le K^2 \ell_{M_t}(\calC)^2.$$ 
\label{theorem:L2bound}
\end{theorem}

To apply Theorems \ref{theorem:mvt} and \ref{theorem:L2bound} to bound
the pointwise norms $\|\eta_t(p)\|$ and $\|*D_t\eta_t(p)\|$ we need to
control the injectivity radius of $p$ and the distance from $p$ to
$\tube^\epsilon_t(\calC)$. 
To bound these two quantities we use the following estimates of
R. Brooks and J. Matelski on the 
geometry of equidistant tori about a short geodesic (see
\cite{Brooks:Matelski:collars}). Their original result only applies to
tubes about non-singular geodesics. The extension to tubes about a
cone singularity is straightforward. 

\begin{theorem}[Brooks-Matelski]
Given $\epsilon \in [0,\varepsilon]$, there are two continuous 
positive functions $d_\epsilon^u$ and
$d_\epsilon^l$ on $[0,\varepsilon]$ 
with $d_\epsilon^u( \delta) \to 0$ as $\delta \to \epsilon$ and
$d_\epsilon^l(\delta) \to
\infty$ as $\delta \to 0$
so that given $\delta \in [0,\varepsilon]$ the
distance between the boundaries of $\tube^\epsilon$ and $\tube^\delta$
satisfies $$d_\epsilon^l(\delta) \le 
d(\bdry \tube^\epsilon,\bdry \tube^\delta)
\le d_\epsilon^u(\delta).$$

\label{theorem:BM}
\end{theorem}

Given Riemannian manifolds $(M,g)$ and $(N,g')$, a diffeomorphism $$h
\colon (M,g) \to (N,g')$$ is {\em $L$-bi-Lipschitz} if we have the
bound $$\sup_{p \in M} {\rm bilip}_p(h^* g', g) \le L.$$ It is worth
noting that we will often be interested in the case when $M=N$ and $h$
is the identity. 

Our $L^2$-bound on $\eta_t$, together with the above
mean value inequality for harmonic strain fields
(Theorem~\ref{theorem:mvt}) and Proposition \ref{curvature}, readily
gives the following corollary. 
\begin{cor} For any $\epsilon>0$, $\delta > 0$, $C>0$ and $L >1$,  
there exists $\ell > 0$ so that if $\ell_{M_\alpha}(\calC) < \ell$ then
the following holds: Let $W$ be a subset of $N$, $\gamma(s)$ a smooth
curve in $W$ and $C(t)$ the geodesic curvature of $\gamma$ in the
$g_t$-metric at $\gamma(0)$. If 
$$W \subset M_t^{\ge \epsilon}$$ for all $t \in [t_0,  \alpha]$ and $C(0)
\le C$ then the identity map
$$\id: (W, g_\alpha) \to{} (W, g_{t_0})$$
is $L$-bi-Lipschitz and $|C(0) - C(A)| \le \delta$. 
\label{cor:bilip}
\end{cor}

To apply the Corollary we need to show that the thick part of
$M_\alpha$ does
not become too thin in $M_t$, while the thin part does not become too thick.
\begin{theorem}
\label{staythick}
Given an $\epsilon_1>0$ there exists
an $\epsilon_0 >0$ and $\ell >0$ so
that if the length $\ell_{M_\alpha}(\calC) < \ell$ then we have 
\begin{equation}
\label{thickone}
M_\alpha^{\ge \epsilon_1} \subset M_t^{\ge \epsilon_0}
\end{equation}
and
\begin{equation}
\label{thicktwo}
M^{\ge \epsilon_1}_t \subset M^{\ge \epsilon_0}_\alpha
\end{equation}
for all $t \in [0,\alpha]$.
\end{theorem}

\bold{Proof:} By Theorem \ref{theorem:BM} we can choose $\epsilon_0>0$
so that $$d(\del M^{\le \epsilon_0}_t, \del M^{\le \epsilon_1/2}_t)
\ge 3\epsilon_1.$$ 
The set $A$ of $t$ such that \eqref{thickone} holds is open in $[0,
\alpha]$. To prove \eqref{thickone} we will show that if
$\ell_{M_\alpha}(\calC)$ is sufficiently short then $A$ is closed. 
Let $a$ be a point in the closure of $A$. By continuity we have
$$M_\alpha^{ \ge 
\epsilon_1} \subseteq M_a^{ \ge \epsilon_0}.$$
Either $a$ is in $A$ and we have proven \eqref{thickone}, or
$M_\alpha^{\ge \epsilon_1} \cap M_a^{\le
\epsilon_0}$ is non-empty.

We work by contradiction and assume that $q \in (M^{\ge
\epsilon_1}_\alpha \cap M^{\le \epsilon_0}_a)$. Let $B$ be a ball of
radius 
$\epsilon_1$ in the $g_\alpha$-metric with $q$ in $\bdry B$ and center
$p$. We also assume that $B$ is contained in 
$M_\alpha^{\ge \epsilon_1}$. By Corollary~\ref{cor:bilip} there exists
an $\ell$ such that if $\ell_{M_\alpha}(\calC) \le \ell$ then the
inclusion map 
$$\iota: (M_\alpha^{\ge \epsilon_1}, g_\alpha) \to{} (M^{\ge
\epsilon_0}_a, g_a)$$ is $2$-bi-Lipschitz. This implies that $p \in
M^{\ge \epsilon_2 /2}_a$ while $d(p,q)$ is less than $2\epsilon_1$ in
the $g_a$-metric.  By our choice of $\epsilon_0$, however, we have $d(p,
\del M^{\le \epsilon_0}_a) \ge 3\epsilon_1$ which contradicts our
assumption that $q$ lies in $M^{\le \epsilon_0}_a$, proving
\eqref{thickone}.

The inclusion \eqref{thicktwo} is proved similarly.
\qed

Before we continue we need to fix some constants. First choose
$\epsilon_2 < \varepsilon$ such that Theorem \ref{theorem:BM} implies
that $d(\del M^{\varepsilon}_t, \del M^{\epsilon_2}_t) >1$. Next
choose $\epsilon_1 < \epsilon_2$ such that $d(\del M^{\epsilon_2}_t,
\del M^{\epsilon_1}_t) > 2 \epsilon_2$. Finally choose $\epsilon_0<
\epsilon_1$ and $\ell_0$ to satisfy the conditions of Theorem
\ref{staythick}. This implies that if $\ell_{M_\alpha}(\calC) \le
\ell_0$ then the inclusion map
$$\iota: (M^{\ge \epsilon_1}_\alpha, g_\alpha) \to{} (N,
g_t)$$ 
is an $L$-bi-Lipschitz diffeomorphism to its image where $L$
only depends on $\ell_{M_\alpha}(\calC)$ and $L \rightarrow 1$ as
$\ell_{M_\alpha}(\calC) \rightarrow 0$. The remainder of this section
will be spent extending this map to all of $M_\alpha \backslash
\tube^{\epsilon_2}_\alpha(\calC)$ in a uniformly bi-Lipschitz way.

Theorem~\ref{theorem:bilip} 
will easily follow from the next result.
\begin{theorem}
\label{intothetube}
Let $V \subset N$ be either
\begin{enumerate}
\item the $\varepsilon$-Margulis tube $\tube^{\varepsilon}_\alpha(\gamma)$
about a geodesic $\gamma$ with $\ell_{M_\alpha}(\gamma) < \epsilon_1$, or
\item a rank-2 cusp component $\cusp^{\varepsilon}_\alpha$ of $M_\alpha^{\le
\varepsilon}$, 
\end{enumerate}
For each $L>1$ there exists and $\ell>0$ such that if
$\ell_{M_\alpha}({\cal C}) \le \ell$ then for all $t \le \alpha$ there
exists an $L$-bi-Lipschitz embedding
$$\phi_t : (V, g_\alpha) \to{} (N, g_t)$$
such that $\phi_t$ is the identity on a neighborhood of $\bdry V$. 
\end{theorem}

We will prove the theorem via a sequence of lemmas.  For simplicity,
these lemmas will treat the case of the Margulis tube; the rank-2 cusp
case admits a far simpler direct proof and is also a limiting case of
these arguments.

We first fix more notation: focusing our attention on a single 
short geodesic $\gamma$, let $W = \tube^{\epsilon_2}_\alpha(\gamma)$
and $T = \del W$.  

\begin{lem}
\label{Tthick}
For each $d>0$ there exists an $\ell>0$ such that if
$\ell_{M_\alpha}({\cal C}) \le \ell$ then $T$ is contained in the
$d$-neighborhood $\calN_d(\del \tube^{\epsilon_2}_t(\gamma))$ of $\del
\tube^{\epsilon_2}_t(\gamma)$ and $W$ contains 
$\tube^{\epsilon_1}_t(\gamma)$ for all $t \le \alpha$. 
\end{lem}

\bold{Proof:} The tubes $\tube_t^{\epsilon_1}(\gamma)$ will vary
continuously in $N$ a $t$ varies. Since $W \supset
\tube_\alpha^{\epsilon_1}(\gamma)$, if $T$ is in $M^{> \epsilon_1}_t$
for all $t$ then $W \supset \tube_t^{\epsilon_1}(\gamma)$ for all 
$t$.  For $d$ sufficiently small $\calN_d(\del
\tube^{\epsilon_2}_t(\gamma))$ will be contained in $M^{>
\epsilon_1}_t$. In particular, the first conclusion implies the
second. 

By Theorem \ref{theorem:BM} there exists an $L>0$ such that if the
injectivity radius in the $g_t$-metric of all points in $T$ lies in
the interval
$[\epsilon_2/L, L\epsilon_2]$ then $T$ is contained in the
$d$-neighborhood of $\del \tube^{\epsilon_2}_t(\gamma)$. 

Let $B$ be a ball of radius $\epsilon_2$ in the $g_\alpha$-metric
centered at a point $p \in T$. Since $B$ is contained in $M^{\geq
\epsilon_1}_\alpha$ 
we have $B \subset M^{\geq \epsilon_0}_t$ for all 
$t$ by Theorem \ref{staythick}. By Corollary \ref{cor:bilip} we can
choose an 
$\ell>0$ such that if $\ell_{M_\alpha} \le \ell$ then the identity map
restricted to $B$ is $L$-bi-Lipschitz from the $g_\alpha$-metric to
the $g_t$-metric.  In particular in the $g_t$-metric there is a ball
of radius $\epsilon_2/L$ centered at $p$ and contained in $B$;
i.e. $p$ has
injectivity radius greater than $\epsilon_2/L$ for all $t$. On the other
hand if $p$ has injectivity radius greater than $L\epsilon_2$ in the
$g_t$-metric then there exists a ball $B'$ of radius greater than
$L\epsilon_2$ in the $g_t$-metric centered at $p$. Reversing the
process above this implies that $p$ has injectivity radius greater
than $\epsilon_2$ in the $g_\alpha$-metric. This contradiction implies
that the injectivity radius at $p$ is bounded above and below by
$L\epsilon_2$ and $\epsilon_2/L$, respectively, for all $t$. \qed

\begin{lem}
\label{stayconvex}
There exists and $\ell > 0$ such that if $\ell_{M_\alpha}({\cal C})
\le \ell$ then $T$ is convex in the $g_t$ metric for all $t \le
\alpha$. 
\end{lem}

\bold{Proof:} To show that $T$ is convex we need to show that every
smooth curve $\sigma$ on $T$ has non-zero geodesic curvature in the
$g_t$-metric at every point on $\sigma$. Since the tube
$\tube_\alpha^{\epsilon_2}(\gamma)$ has radius uniformly bounded below
by some $R>0$, every smooth curve on $T$ has geodesic curvature
greater than $\tanh R$ in the $g_\alpha$-metric.  The result then
follows from Corollary \ref{cor:bilip}. \qed

We denote by 
$$\pi_t: T \to{} \tube^{\epsilon_2}_t(\gamma)$$ the radial
projection mapping. More explicitly, each point $p \in T$ lies on a
geodesic ray which is perpendicular to the core of
$\tube^\varepsilon_t(\gamma)$ in the 
$g_t$-metric. This ray intersects $\del \tube^{\epsilon_2}_t(\gamma)$ in a
unique point $p'$, and we set $\pi_t(p) = p'$. 

\begin{lem}
For each $L>1$ and $\delta>0$ there exists an $\ell>0$ such that if
$\ell_{M_\alpha}({\cal C}) \le \ell$ then 
\begin{enumerate}
\item the radial projection
$\pi_t$ is an
$L$-bi-Lipschitz diffeomorphism, and
\item if $\sigma$ is a geodesic
in the Euclidean metric on $T$ induced by the $g_\alpha$-metric then
$\pi_t(\sigma)$ has curvature bounded by $\delta$ in the Euclidean
metric on $\del \tube^{\epsilon_2}_t(\gamma)$ induced by the $g_t$-metric. 
\end{enumerate}
\label{lemma:geodesic:curvature}
\end{lem}

\bold{Proof:} Given $p \in T$, let $P$ be the hyperbolic plane in the
$g_t$-metric tangent to $T$ at $p$, and let $\ray$ be the radial geodesic
through $p$. By Lemma \ref{Tthick} a bound $\ell_{M_\alpha}(\calC)$
gives a bound on $d(p, \pi_t(p))$ and as $\ell_{M_\alpha}(\calC)
\rightarrow 0$ we have $d(p, \pi_t(p)) \rightarrow 0$. We show
that a bound on $\ell_{M_\alpha}(\calC)$ gives a bound on the angle
$\angle(\ray, P)$ between $\ray$ and $P$ and as
$\ell_{M_\alpha}(\calC) \rightarrow 0$ we have $|\angle(\ray, P) -
\pi/2| \rightarrow 0$. 

To control $\angle(\ray, P)$ we make the following
observation. Let $$W'=\tube_t^{\epsilon_2}(\gamma) \backslash
\calN_d(\bdry \tube_t^{\epsilon_2}(\gamma)).$$ 
By Lemma \ref{Tthick} we have $W \supset W'$ since $T$ is convex $P$
and $W'$ are disjoint. On the other hand $p \in P$ is within $2d$ of
$W'$ and the tube $W'$ has definite radius. Elementary hyperbolic
geometry then gives the desired bound.

Next we remark that if $\angle(\ray, P) \neq 0$ then $\pi_t$ is a
diffeomorphism at $p$. If $\pi_t$ is a local diffeomorphism at each
point in $T$ it is a covering map. Since $T$ and $\del 
\tube_t^{\epsilon_2}(\gamma)$ are homotopic in the complement of the
core geodesic, $\pi_t$ must be a global diffeomorphism. To finish the
proof of (1) we note that bounds on $d(p, \pi_t(p))$ and
$|\angle(\ray, P) - \pi/2|$ along with a lower bound on the tube
radius of $W'$ bound the bi-Lipschitz constant of $\pi_t$ at $p$. 
 
Next we control the curvature of the curve $\overline{\sigma} = \pi_t 
\circ \sigma$. We give
the tube $\tube_t^\varepsilon(\gamma)$ cylindrical coordinates $(r,
\theta, z)$ so that $g_t\vert_{\tube_t^\varepsilon(\gamma)}$ is given by
the Riemannian metric $$dr^2 + \sinh^2 r d\theta^2 +
\cosh^2r dz^2,$$ 
where $r$ measures the hyperbolic distance from the core geodesic
of $\tube_t^\varepsilon(\gamma)$. We then let $\sigma(s) = (r(s),
\theta(s), z(s))$ be a unit speed parameterization of $\sigma$.

We begin the proof of (2) with some preliminary remarks. 

The bound on $|\angle(\ray, P) - \pi/2|$ described above gives a bound
on $r'(s)$.  We also note that by the remark after
Proposition~\ref{curvature}, $\sigma$ will be almost geodesic on $(T,
g_t)$. In particular $\nabla^t_{\sigma'} \sigma'$ will be nearly
orthogonal to $T$ and hence nearly radial. 
That is, we can bound
$\angle(\nabla^t_{\sigma'} \sigma'(s), \frac{\del}{\del r})$. Putting
this all together, for any $\epsilon>0$ if $\ell_{M_\alpha}(\calC)$ is
sufficiently small, then $|r'(s)| \le \epsilon$ and 
$$\angle(\nabla^t_{\sigma'}\sigma'(s), \frac{\del}{\del r}) \le
\epsilon.$$ 

If $\tube_\alpha^{\epsilon_2}(\gamma)$ has radius $R_\alpha$ then the
curvature of $\sigma$ in the $g_\alpha$-metric will be less than
$\coth R_\alpha$. By Corollary \ref{cor:bilip} for any $C>0$ if
$\ell_{M_\alpha}(\calC)$ is sufficiently small then the curvature of
$\sigma$ in the $g_t$-metric will be less than $\coth R_\alpha +
C$.  The curvature is the length of $\nabla^t_{\sigma'}
\sigma'$ in the $g_t$-metric, so we have
$$\left| \nabla^t_{\sigma'} \sigma' \right|_t \le \coth R_\alpha +
C.$$ Combining this with our bound $\angle(\nabla^t_{\sigma'}\sigma',
\frac{\del }{\del R}) \le \epsilon$ 
we have
$$\left|g_t\left(\nabla^t_{\sigma'}\sigma', \frac{1}{\sinh r}
\frac{\del}{\del \theta}\right)\right| \le \epsilon (\coth R_\alpha +
C).$$ 

Direction calculation gives
\begin{eqnarray*}
\nabla^t_{\sigma'}{\sigma'}(s) & = & \left(r''(s) - \sinh r(s)
\cosh r(s) \left( \theta'(s)^2 +
z'(s)^2 \right) \right)\frac{\del}{\del r} \\ 
& & + \left( \theta''(s) + 2\coth r(s)r'(s) \theta'(s) \right)
\frac{\del}{\del \theta} \\ 
& & + \left(z''(s) + 2\tanh r(s) r'(s) z'(s) \right) \frac{\del}{\del z}
\end{eqnarray*}
so
$$\left|g_t\left(\nabla^t_{\sigma'} \sigma', \frac{1}{\sinh r}
\frac{\del }{\del \theta}\right)\right| = \left|\sinh r \theta'' + 2 \cosh
r r' \theta'\right| \le \epsilon (\coth R_\alpha + C).$$ Since 
$\sigma$ is unit speed 
we have $|\cosh r \theta' | \le |\coth r|$.  Combining this with the
fact that $|r'| < \epsilon$ we have
$$|\sinh r \theta''| \le \epsilon  (2\coth r + \coth R_\alpha + C).$$
By a similar method we also see that
$$|\cosh r z''| \le \epsilon  (2 \tanh r + \coth R_\alpha + C).$$

We can now bound the curvature of $\overline{\sigma}$ on $(\del
\tube_t^{\epsilon_2}(\gamma), g_t)$. We first note that 
if $\overline{\nabla}^t$ is the Riemannian connection for the
Euclidean metric on $\del \tube_t^{\epsilon_2}(\gamma)$ then
$$\overline{\nabla}^t_{\overline{\sigma}'} \overline{\sigma}'(s) =
\theta''(s) \frac{\del}{\del \theta} + z''(s) \frac{\del}{\del z}.$$
Since $\overline{\sigma}$ is not necessarily unit speed, the length of
$\overline{\nabla}^t_{\overline{\sigma}'} \overline{\sigma}'(s)$ does
not necessarily give the geodesic curvature of $\overline{\sigma}$.
If, however, $\overline{\sigma}(h(s))$ is a unit speed
reparameterization of $\overline{\sigma}$ then $h'(s)$ is close to 1
and $h''(s)$ is small.

Thus, the geodesic curvature of 
$\overline{\sigma}$ on $(\del \tube_t^{\epsilon_2}(\gamma), g_t)$ 
is approximately the norm of
$\overline{\nabla}^t_{\overline{\sigma}'} \overline{\sigma}'(s)$
which is bounded by
$$|\overline{\nabla}^t_{\overline{\sigma}'} \overline{\sigma}'(s) |^2
= \sinh^2 r (\theta'')^2 + \cosh^2 (z'')^2 \le 2\epsilon^2(2\coth
r  + \coth R_\alpha + C)^2.$$  The right hand side limits to
zero as $\ell_{M_\alpha}(\calC) \rightarrow 0$, completing the proof.
\qed

We now return to the situation at hand, and recall that the inclusion
$$\iota \colon (M_\alpha^{\ge \epsilon_1}, g_\alpha) \to (N,g_t)$$ 
is an
$L$-bi-Lipschitz diffeomorphism to its image.  Given the short
geodesic $\gamma$ with length 
$$\ell_{M_\alpha}(\gamma) \le \epsilon_1,$$ we now show 
that the inclusion map can be modified to a map $\phi_t$ on the collar
$\tube^\varepsilon_\alpha(\gamma) \setminus
\tube^{\epsilon_2}_\alpha(\gamma)$ that is bi-Lipschitz and sends $\bdry
\tube^{\epsilon_2}_\alpha(\gamma)$ to $\bdry \tube^{\epsilon_2}_t(\gamma)$.
\begin{lem}
For each $L>1$ there exists an $\ell>0$ such that if
$\ell_{M_\alpha}({\cal C}) \le \ell$ then there exists an embedding
$$\phi_t: 
(\tube^\varepsilon_\alpha(\gamma) \backslash
\tube^{\epsilon_2}_\alpha(\gamma), g_\alpha) 
 \to{} (N, g_t)$$ 
 such that 
\begin{enumerate}
\item $\phi_t$ is the identity in a neighborhood of $\del
\tube_\alpha^\varepsilon(\gamma)$, 

\item $\phi_t$ is $L$-bi-Lipschitz from the $g_\alpha$-metric to the
$g_t$-metric, and

\item $\phi_t = \pi_t$ on $T$.

\item $d \phi_t$ sends each unit normal vector to
$\bdry \tube^{\epsilon_2}_\alpha(\gamma)$ to a unit normal vector
to $\bdry \tube_t^{\epsilon_2}(\gamma)$.
\end{enumerate} 
\label{lemma:tube:collar} 
\end{lem} 

\bold{Proof:}  Using a smooth bump function, it is straightforward to
extend the projection $\pi_t$ to a map 
$$\phi_t(p,r) = (p,s(r))$$ 
where if $(p,r)$ lies in $T$ then $(p,s(r)) = \pi_t(p)$ and
so that $\phi_t$ is the identity on $\bdry
\tube^{\epsilon_2}_\alpha(\gamma)$.  
Since by Lemma~\ref{Tthick}, for any $d>0$ we may choose $\ell>0$ 
so that the Margulis tube
$\tube^{\epsilon_2}_t(\gamma)$ has boundary
$\bdry \tube^{\epsilon_2}_t(\gamma)$ lying within distance $d$ of $T =
\bdry \tube^{\epsilon_2}_\alpha(\gamma)$ in $N$, 
$s(r)$ may be chosen so that $\phi_t$ is $L$-bi-Lipschitz.  

Since radial geodesics in $\tube_\alpha^{\varepsilon}(\gamma)$ make
small angle with $\bdry \tube_t^{\epsilon_2}(\gamma)$
(as in the proof of Lemma~\ref{lemma:geodesic:curvature}), a further small
modification in a neighborhood of $T$ ensures that $d\phi_t$ sends
unit normal vectors to $\bdry \tube^{\epsilon_2}_\alpha(\gamma)$ to
unit normal vectors to $\bdry \tube^{\epsilon_2}_t(\gamma)$.
\qed

The shape of the Margulis tube $\tube^{\epsilon}(\gamma)$ about a
short geodesic $\gamma$ in a hyperbolic 3-manifold varies continuously
with the complex length of the core curve $\gamma$ in
the smooth bi-Lipschitz topology (cf. \cite[Lem. 6.2]{Minsky:torus}). By
\cite[Thm. 4.3]{Bromberg:Schwarzian}, when the cone-singularity is
sufficiently short one can control the derivative of the complex
length of a bounded length closed geodesic in $M_t$ (for the analogous
statement for the Teichm\"uller parameter for a rank-2 cusp, see
\cite[Thm. 7.3]{Bromberg:Schwarzian}).  Together, these 
estimates yield the following lemma:
\begin{lem}
For each $L>0$ there exists an $\ell>0$ such that if
$\ell_{M_\alpha}({\cal C}) \le \ell$ then there exists an
$L$-bi-Lipschitz diffeomorphism 
$$\psi_t :\tube^{\epsilon_2}_\alpha(\gamma) \to{}
\tube^{\epsilon_2}_t(\gamma)$$
so that $\psi_t$ restricts to $\bdry \tube_\alpha^{\epsilon_2}(\gamma)$ by
an affine map with respect to the Euclidean structures on 
$\bdry \tube_\alpha^{\epsilon_2}(\gamma)$ and
$\bdry \tube_t^{\epsilon_2}(\gamma)$ and $d\psi_t$
sends unit normal vectors to 
$\bdry \tube_\alpha^{\epsilon_2}(\gamma)$ to unit normal vectors to
$\bdry \tube_t^{\epsilon_2}(\gamma)$.
\label{lemma:continuous:tubes}
\end{lem}

Applying Lemma~\ref{lemma:tube:collar}, the final step in
our argument will be the following lemma.

\begin{lem}
Assume $\tube^{\epsilon_2}$ is an $\epsilon_2$-Margulis tube whose
core geodesic has length less than $\epsilon_0$. Let  
$$\phi \colon \bdry \tube^{\epsilon_2} \to \bdry \tube^{\epsilon_2}$$
be a diffeomorphism isotopic to the identity.  
Then for each $K>1$ 
there is an $L>1$ and $\delta>0$ 
so that if $\phi$ is $L$-bi-Lipschitz and sends geodesics in the
intrinsic metric on 
$\bdry \tube^{\epsilon_2}$ to arcs of geodesic curvature at
most $\delta$ in the intrinsic metric on $\bdry
\tube^{\epsilon_2}$, 
then $\phi$ extends to a $K$-bi-Lipschitz diffeomorphism $$\Phi \colon
\tube^{\epsilon_2} \to \tube^{\epsilon_2}.$$
\label{lemma:tube:extension}
\end{lem}

\bold{Proof:}   
If $\tube^{\epsilon_2}$ were instead a rank-2 cusp component
$\cusp^{\epsilon_2}$ of $M^{\le \epsilon_2}$,
it would admit a natural
parameterization $T \times \reals^+$ where $(x,d) \in T \times
\reals^+$ represents the point at depth $d$ along inward-pointing
normal to $T$ at $x \in T$.  Since the radial projection $\pi_d \colon
T \to T_d$ from $T$ to the torus $T_d$ at depth $d$ given by
$$\pi_d(x,0) = (x,d)$$ is conformal, the radial extension $$\Phi(x,d)
= (\phi(x),d)$$ is readily seen to be bi-Lipschitz, with bi-Lipschitz
constant $L$.

For the Margulis tube $\tube^{\epsilon_2}$ the situation is
very similar away from a neighborhood of the core geodesic.  Indeed,
after removing the core geodesic, the tube has a natural product
structure $T \times \reals^+$, 
where $(x,r) \in T \times
\reals^+$  is now a point at radius $r$ from $\gamma$.  If
$\tube^{\epsilon_2}$ has radius $R$, the radial
projection
$$\pi_r(x,R) = (x,r),$$
for $r \le R$, 
is uniformly {\em quasi-conformal} for $r \ge 1$,
and we will see that the radial
extension 
$$\Phi(x,r) =
(\phi(x),r)$$
is uniformly bi-Lipschitz for $r \ge 1$ 
as a result.  For $r \in (0,1)$, however, the radial extension 
$\Phi$ on over $\tube^{\epsilon_2}$ will not in general be
bi-Lipschitz, nor will it in general extend to the core geodesic.

To correct these problems, we will construct an isotopy $\phi_r$ of
$\phi$ to the identity and define our extension $\Phi$ by
\begin{equation}
\Phi(x,r) = (\phi_r(x),r).
\label{extension}
\end{equation}
The isotopy $\phi_r$ will be the flow
of a time dependent family of vector fields $v_r$ on $T$.  If $T_r$
denotes the equidistant torus at radius $r$ from $\gamma$, the
extension $\Phi$ will be bi-Lipschitz if each $\phi_r$ is uniformly
bi-Lipschitz on $T_r$ and the vector field $v_r$ has uniformly bounded
size on $T_r$.

Using the intrinsic Euclidean metric on $T_r$, we identify all tangent
spaces $T_x(T_r)$, $x \in T_r$, with $\reals^2$ via parallel
translation. In particular for each $x$ we view $d\phi_x$ as a linear
map of $\reals^2$ to itself. We choose an orthonormal framing $\{e_1,
e_2\}$ such that $e_1$ and $e_2$ are tangent to the directions of
principal curvature of $T_r$. Then $| d\phi - \id|_r$ is the maximum
of the absolute value of the entries of the matrix $d\phi - \id$
written it terms of this basis. 

Let $\tube^{\epsilon_2}$ have radius $R$ so that $T_R = \bdry
\tube^{\epsilon_2}$. 
We first claim that given any
$\delta'>0$ there are $\delta > 0$ and 
$L>1$ such that $|d\phi - \id|_R \le \delta'$. This is equivalent to
showing that $\phi$ has small bi-Lipschitz constant and small
``twisting.''  That is we need to show that the angle between any $v$
and $d\phi(v)$ is small for any tangent vector $v$. We have assumed
that $\phi$ is $L$-bi-Lipschitz so we only need to bound the twisting.

Let $\alpha$ be the shortest geodesic on $T_R$ through $x$.
Then length of $\alpha$ is bounded by $C\epsilon_2$, where $C$ is a
universal constant, so the length of $\phi(\alpha)$ is bounded by
$L C \epsilon_2$. By our assumption the geodesic curvature of $\phi(\alpha)$ is
bounded by $\delta$.  Since $\phi(\alpha)$ is homotopic to $\alpha$, 
there is some point $y$ on $\alpha$ such that the if $v$ is tangent to
$\alpha$ at $y$ then the angle between $v$ and $d\phi(v)$ is $0$. 
The bounds on the curvature and the length of $\phi(\alpha)$ imply
that the tangent to $\phi(\alpha)$ is nearly parallel to $\alpha$
everywhere. Therefore, for any vector $v$ tangent to $\alpha$, $v$ and
$d\phi(v)$ make a small angle. Since $L$ is close to $1$,  $\phi$ is
nearly conformal, so for any tangent vector $v$ the angle between $v$
and $d\phi(v)$ is small, so $\phi$ has small twisting.

We now consider a linear homotopy of $\phi$ to the identity
constructed as follows. After normalizing by an isometry of
$\tube^{\epsilon_2}$, we may assume that $\phi$ fixes a point $p
\in T_R$. We then identify the universal cover $\wt T_R$
with $\reals^2$ so that the intrinsic metric on $T_R$ lifts to the
Euclidean metric on $\reals^2$.  We let $\wt \phi$ be a 
lift of $\phi$ that fixes a point and hence a lattice.
Let $\wt \phi_t$ be the homotopy of $\wt \phi$ to the identity given by 
$$\wt \phi_t (\vec{x}) = (1 - t) \wt \phi(\vec{x}) + t \vec{x}$$ 
for each $\vec{x}$ in $\reals^2$. This homotopy is equivariant by 
the action of the covering translation group for $T_R$ and therefore
descends to a homotopy $\phi_t$ of $T_R$. Direct computation shows
that if $|d\phi - \id|_R \leq \delta'$ then $|d\phi_t - \id|_R \leq
(1-t)\delta'$. In particular, for $\delta'$ sufficiently small,
$\phi_t$ is a local, and hence global, diffeomorphism for all $t$. In
the 
$T_r$ metric another computation shows that 
$$|d\phi - \id|_r \le \frac{\tanh R}{\tanh r}\delta'.$$ Therefore for
any $\delta'' \ge 0$ we can chose $L$ and $\delta$ such that $|d\phi_t
- \id| \le \delta''$ in the $T_r$-metric for all $r\ge 1$.  A bound on
$|d\phi_t - \id|_r$ determines a bound on the bi-Lipschitz
constant. In particular, for any $L'>1$ we can choose $L$ and $\delta$
so that 
$\phi_t$ is $L'$-bi-Lipschitz in the $T_r$-metric.

Let $v_t$ be the time dependent family of vector fields whose flow is
$\phi_t$. The norm of $v_t$ in the $T_r$ metric is bounded by the
suprememum of the distance between $x$ and $\wt \phi(x)$ in the $T_r$
metric.  If $D_r$ is the diameter of $T_r$, this distance is bounded
by $D_r \epsilon'$ where $\epsilon'$ tends to zero as the bi-Lipschitz
constant $L'$ tends to 1.  But $D_r$ is universally bounded for $1 \le
r \le 2$ so for these values of $r$ the norm of $v_t$ is bounded on
$T_r$.

To finish the proof, we must reparameterize our isotopy to obtain the
desired smoothness for the extension $\Phi$. Let $s$ be
a smooth function on $[0,R]$ with
$$s(r) = \left\{ \begin{array}{ll}
1 & r \le 1 \\
\mbox{monotonically decreasing} & 1 \le r \le 2 \\
0 & 2 \le  r \le R
\end{array}
\right.
$$
Note that the derivative of $s$ is bounded (and independent of $R$).
We now abuse notation and redefine the isotopy $\phi_r$ as the
projection to $T$ of the isotopy of $\wt \phi$ to the identity 
given by the formula 
$$\wt \phi_r(\vec{x}) = (1-s(r))\wt \phi(\vec{x}) + s(r)\vec{x}$$
for $r \in [0,R]$.
With this notation,
$\phi_r$ is the flow of the time dependent vector field 
$$w_r = s'(r) v_{s(r)}$$ which is again bounded in the $T_r$-metric
for all $r$, and zero for $r \le 1$ and $r \ge 2.$ It follows that the
extension $\Phi$ of $\phi$ given by (\ref{extension}) is
$K$-bi-Lipschitz, where $K$ depends only on $L$ and $\delta$.  The
proof is complete.\qed

We now prove Theorem~\ref{intothetube} by combining
Lemma~\ref{lemma:tube:collar}, 
Lemma~\ref{lemma:continuous:tubes}, and
Lemma~\ref{lemma:tube:extension}. 

\bold{Proof:} Let $L >1$ and $\delta >0$ be given.  Then by 
Lemma~\ref{lemma:tube:collar} there is an $\ell>0$
so that if $\ell_{M_\alpha}(\calC) \le \ell$ then for each closed
geodesic $\gamma$ with $\ell_{M_\alpha}(\gamma) < \epsilon_1$ we have
the
corresponding
mapping
$$\phi_t: 
(\tube^\varepsilon_\alpha(\gamma) \backslash
\tube^{\epsilon_2}_\alpha(\gamma), g_\alpha) 
 \to{} (N, g_t).$$ 

Then $\phi_t$ restricts to an $L$-bi-Lipschitz diffeomorphism  
$$\phi_t\vert_{\bdry \tube_\alpha^{\epsilon_2}(\gamma)} \colon
\bdry \tube_\alpha^{\epsilon_2}(\gamma) \to
\bdry \tube_t^{\epsilon_2}(\gamma)$$
so that geodesics on $\bdry \tube_\alpha^{\epsilon_2}(\gamma)$ map to
arcs of geodesic curvature bounded by $\delta$ on 
$\bdry \tube_t^{\epsilon_2}(\gamma)$.

Assume $\ell$ also satisfies the hypotheses of
Lemma~\ref{lemma:continuous:tubes} and let
$$\psi_t \colon \tube_\alpha^{\epsilon_2}(\gamma) \to
\tube_t^{\epsilon_2}(\gamma)$$
be the $L$-bi-Lipschitz diffeomorphism guaranteed by
Lemma~\ref{lemma:continuous:tubes}.

Since 
$$\psi_t \colon \bdry \tube_\alpha^{\epsilon_2}(\gamma) \to
\bdry \tube_t^{\epsilon_2}(\gamma)$$ is affine, 
the composition
$$\phi_t \compos \psi_t^{-1} \vert_{\bdry
\tube_t^{\epsilon_2}(\gamma)} \colon
\bdry \tube_t^{\epsilon_2}(\gamma) \to
\bdry \tube_t^{\epsilon_2}(\gamma)$$
sends geodesics on $\bdry \tube_t^{\epsilon_2}(\gamma)$
to arcs of curvature bounded by $\delta$, and is $L^2$-bi-Lipschitz.

Since $L>1$ and $\delta$ were arbitrary, given any $K>1$ we may choose $L$ and
$\delta$ sufficiently small so that 
Lemma~\ref{lemma:tube:extension} provides a $K$-bi-Lipschitz extension
$$\Phi_t \colon \tube_t^{\epsilon_2}(\gamma) \to
\tube_t^{\epsilon_2}(\gamma)$$ of 
$\phi_t \compos \psi_t^{-1} \vert_{\bdry
\tube_t^{\epsilon_2}(\gamma)}$
over $\tube_t^{\epsilon_2}(\gamma)$.

Since we 
have 
$$\phi_t \vert_{\bdry \tube_\alpha^{\epsilon_2}(\gamma)} =
\Phi_t^{-1} \compos \psi_t \vert_{\bdry \tube_\alpha^{\epsilon_2}(\gamma)}$$
the composition
$$\Phi_t^{-1} \compos \psi_t \colon \tube_\alpha^{\epsilon_2}(\gamma)
\to \tube_t^{\epsilon_2}(\gamma)$$ gives a 
$KL$-bi-Lipschitz extension of $\phi_t$ over
$\tube_\alpha^{\epsilon_2}(\gamma)$.

As remarked, we may apply exactly analogous versions of 
Lemma~\ref{lemma:tube:collar}, 
Lemma~\ref{lemma:continuous:tubes}, and
Lemma~\ref{lemma:tube:extension} for the rank-2 cusp case to complete
the proof.
\qed

\bold{Proof:} {\em (of Theorem~\ref{theorem:bilip})}.
Given any $L>1$, we may choose $\ell>0$ satisfying the hypotheses of
Theorem~\ref{intothetube} so that
the inclusion
$$\iota \colon (M_\alpha^{\ge \epsilon_1}, g_\alpha) \to (N,g_t)$$
is an $L$-bi-Lipschitz diffeomorphism to its image from the
$g_\alpha$-metric to the $g_t$-metric.  After decreasing $\ell$ if
necessary, 
Theorem~\ref{theorem:bilip} follows by applying
Theorem~\ref{intothetube} to extend $\iota$ to an $L$-bi-Lipschitz
diffeomorphism $f_t$ by extending over each component of $M_\alpha^{\le
\varepsilon}$ other than $\tube^{\varepsilon}(\calC)$.  

Finally, for each component $c \subset \calC$ we may apply the
argument of Lemma~\ref{lemma:tube:collar} to each component of 
$\tube^{\varepsilon}(\calC)$ to modify the restriction
$$f_t \vert_{M_\alpha \setminus \tube_\alpha^{\epsilon_2}(\calC)}$$ on
$\tube_\alpha^{\varepsilon}(c) \setminus \tube_\alpha^{\epsilon_2}(c)$
so that 
$f_t (\bdry  \tube_\alpha^{\epsilon_2}(c)) = 
\bdry  \tube_t^{\epsilon_2}(c)$.
The resulting $L$-bi-Lipschitz diffeomorphism 
$$h_t \colon M_\alpha \setminus \tube_\alpha^{\epsilon_2}(\calC)  \to 
M_t \setminus \tube_t^{\epsilon_2}(\calC)$$
$h_t$ proves
Theorem~\ref{theorem:bilip}. 
\qed

We remark that the choice of $\epsilon_2$ was arbitrary, and there
exists $\ell>0$ satisfying the theorem for any choice of $\epsilon_2$.

\bold{Constants.} For future reference, choosing some $L>1$ and
letting $\ell >0$ denote the corresponding constant
so that the conclusions of Theorem~\ref{theorem:bilip} apply, we 
take as a threshold constant
$$\ell_0 = \min \{\ell_{\rm knot}, \ell\}.$$
Then for any $M \in AH(S)$ and any geodesic $\gamma \subset M$ with
$\ell_M(\gamma) < \ell_0$, we may graft $M$ along $\gamma$, and we
may decrease the cone-angle along $\gamma$ with $L$-bi-Lipschitz
distortion of the metric outside of a standard tube about $\gamma$.

\section{Realizing ends on a Bers boundary}
\label{section:realized}
In this section we define a notion of {\em realizability} for
simply degenerate ends of hyperbolic 3-manifolds 
on Bers boundaries for Teichm\"uller space: 
\begin{defn}
Let $E$ be a simply degenerate end of a hyperbolic 3-manifold
$M$.  If $E$ admits a marking and orientation preserving bi-Lipschitz
diffeomorphism to an end $E'$ of a manifold $Q$ lying in the boundary
of a Bers slice, we say $E$ is 
{\em realized on a Bers boundary by $Q$}.
\end{defn}
We prove a general realizability result for ends of
manifolds $M \in AH(S)$.
\begin{theorem}{\sc Ends are Realizable}
Let $M \in AH(S)$ have no cusps.  Then each degenerate end $E$ of $M$
is realized on a Bers boundary by a manifold $Q$.  
\label{theorem:ends:realized}
\end{theorem}

\bold{Proof:}
As usual, the proof breaks into cases.

\ital{Case I: the end $E$ has arbitrarily short geodesics.}
Let $\{\gamma_n\}_{n=0}^\infty$ be a collection of arbitrarily short
geodesics in $M$.  Pass to a subsequence so that each $\gamma_n$ has
length less than $\ell_0$, and so that for each $n \ge 1$ the
geodesic $\gamma_n$ is isotopic out the end $E$ in the complement of
$\gamma_0$.

Consider the simultaneous graftings
$$M^c_n = {\rm Gr}^\pm(\gamma_0,\gamma_n,M).$$
The manifolds $M^c_n$ are 3-dimensional hyperbolic cone-manifolds with
cone-angles $4 \pi$ at the cone-singularities $\gamma_0$ and
$\gamma_n$.

Since we have $$\ell_M(\gamma_n) < \ell_0$$ for each $n$, we may
apply Theorem~\ref{theorem:bilip} to decrease the cone-angles at
$\gamma_0$  and $\gamma_n$ to $2\pi$.
The result is a smooth, geometrically finite hyperbolic 3-manifold
homotopy equivalent to $S$, namely, a quasi-Fuchsian manifold $Q$.  We let
$X_n$ and $Y_n$ in $\Teich(S)$ be the surfaces simultaneously
uniformized by $Q$, so that 
$Q = Q(X_n,Y_n)$.

We note that the
conformal boundary component that arises from negative grafting along
$\gamma_0$ does not change with $n$; there is a
single $X \in \Teich(S)$ so that $X_n = X$.  
It follows that the 
sending the cone-angles of $M^c_n$ to $2 \pi$ gives the sequence
$\{ Q(X,Y_n) \}_{n=1}^\infty$ lying in the Bers slice $B_X.$

Passing to a subsequence and extracting a limit $Q\in \bdry B_X$,
we  claim that the manifold $Q$ realizes the end $E$ on the Bers
boundary $\bdry B_X$. 

To see this, let $U_n$ be the union of Margulis tubes
$$U_n  = \tube^{\varepsilon}(\gamma_0) \disjunion
\tube^{\varepsilon}(\gamma_n)$$ in $M$.
Let
$$F\colon S \times \reals \to M$$ be a
product structure for $M$ as in section~\ref{section:grafting} so that
$\gamma_n$ is a simple curve on $F(S\times \{n\})$.
We consider an exhaustion of the end $E$ by compact submanifolds $K_n
= F (S \times [t_0,t_n])$ where $[t_0,t_n]\subset [0,\infty)$ is an
interval chosen so that $$K_n \cap U_n = \nullset.$$ Letting $F
 (\gamma_0 \times (-\infty,0])$ be the negative grafting annulus for
$\gamma_0$ and $F( \gamma_n \times [n ,+\infty))$ be the positive
grafting annulus for $\gamma_n$,  $K_n$ admits a marking
preserving isometric embedding $\iota_n$ to the subset $\iota_n(K_n) =
K_n' \subset M^c_n$.

Let $U_n' \subset Q(X,Y_n)$ denote the
union $$U_n' = \tube^{\varepsilon}(\gamma_0) \disjunion
\tube^{\varepsilon}(\gamma_n),$$ 
and let $$h_n \colon \left( 
M^c_n \setminus U_n , \bdry U_n \right) 
\to \left( Q(X,Y_n) \setminus U_n', \bdry U_n' \right)$$
be the uniformly bi-Lipschitz diffeomorphisms furnished by
Theorem~\ref{theorem:bilip}.  For each integer $j$, there is an $n_j$
so that the mappings
$\varphi_n = h_n \compos \iota_j$ are uniformly bi-Lipschitz
embeddings of $K_j$ 
into $Q(X,Y_n)$ for all $n > n_j$.  

By Ascoli's theorem, the embeddings $\varphi_n$ converge to a
uniformly bi-Lipschitz embedding on $K_j$ into the limit $Q$ of
$Q(X,Y_n)$ after passing to a subsequence.  Diagonalizing, we have a
uniformly bi-Lipschitz embedding of $E$ into $Q$, so $E$ is realized
by $Q$ on the Bers boundary $B_X$.

\ital{Case II: the end $E$ has bounded geometry.}
By our application of Minsky's bounded geometry results
(Theorem~\ref{theorem:bounded:geometry}) if $Q \in \bdry B_X$ has
end-invariant $\nu(E_Q) = \nu(E)$, then $Q$ realizes the end $E$ on
the Bers boundary $B_X$.  Letting $\mu$ be any measure on $\nu(Q)$ and
choosing weighted simple closed curves $t_n \gamma_n$ so that $t_n
\gamma_n \to \mu$, we let $Y_n \in \Teich(S)$ be any surface for which
$\ell_{Y_n}(\gamma_n) < 1/n$.  Then applying
\cite[Thm. 2]{Brock:length},  any accumulation point $Q$ of $Q(X,Y_n)$
has the property that $\nu(E_Q) = \nu(E)$. 
This completes the proof.
\qed

Applying results of \cite{Brock:length} and \cite{Brock:boundaries}, we
deduce a corollary that will be essential to ensure convergence of our
candidate geometrically finite approximates for $M$.
\begin{cor}
Let $Q = \lim Q(X,Y_n)$ be a realization of the positive degenerate
end $E$ of $M \in AH(S)$.  Then for any convergent subsequence $Y_n
\to [\mu]$ in Thurston's compactification $\pl(S)$, we have $\nu(E) =
|\mu|$.
\label{corollary:pml}
\end{cor}
\bold{Proof:}  The corollary is a direct consequence of Theorem~6.1 of 
\cite{Brock:boundaries}. \qed

\section{Asymptotic isolation of ends}
\label{section:isolation}
Theorem~\ref{theorem:ends:realized} 
guarantees that each degenerate end of a manifold $M \in AH(S)$ can be
realized on a Bers boundary.  In the doubly degenerate case, 
realizations of $E^-$ and $E^+$ by $$ Q^- = \lim_n Q(X_n,Y)
\ \ \ \text{and} \ \ \ 
Q^+ = \lim_n Q(X,Y_n)
$$
suggest the candidate approximates $Q(X_n,Y_n)$ for the original
manifold $M$.  

The construction raises the natural question: to what extent do the
ends of a limit $N$ of quasi-Fuchsian manifolds $Q(X_n,Y_n)$ 
depend on the {\em pair} of surfaces $(X_n,Y_n)$?  In this section we
isolate the effect of the negative surfaces $X_n$ on the positive end
of $N$ and likewise for the negative end.

\begin{theorem}{\sc Asymptotic Isolation of Ends}
Let $Q(X_n,Y_n) \in AH(S)$ be a sequence of quasi-Fuchsian manifolds
converging algebraically to the cusp-free limit manifold $N$.  Then,
up to marking and orientation preserving bi-Lipschitz diffeomorphism,
the positive end of $N$ depends only on the sequence $\{ Y_n \}$ and
the negative end of $N$ depends only on the sequence $\{ X_n \}$.
\label{theorem:isolation}
\end{theorem}

\bold{Proof:}
Consider such a sequence $Q(X_n,Y_n)$ converging to $N$.
By the assumption that $N$ is cusp-free, the convergence of
$Q(X_n,Y_n)$ to $N$ is strong (see \cite{Thurston:book:GTTM} or
\cite[Cor. G]{Anderson:Canary:cores}).   

We show that the positive end $E^+$
depends only on the sequence $\{Y_n\}$ up to bi-Lipschitz
diffeomorphism preserving orientation and marking; the same argument
applies to the negative end $E^-$.  In other words, given another 
sequence $\{ X_n' \} $ in $\Teich(S)$ for which the quasi-Fuchsian manifolds
$Q(X_n', Y_n)$ converge to $N'$, the positive end $(E^+)'$ of $N'$
admits a marking and orientation preserving bi-Lipschitz
diffeomorphism with $E^+$.

If the end $E^+$ is geometrically finite, it is well known that its
associated conformal boundary component $Y$ determines $E^+$ up to
bi-Lipschitz diffeomorphism (see \cite{Epstein:Marden:convex},
\cite{Minsky:ends}).  Thus we may assume that $E^+$ is degenerate.

If all closed geodesics in $N$ have length at least $\ell_0/2$
then $N$ has a global lower bound to its injectivity radius since $N$
has no cusps.  
It follows from Theorem~\ref{theorem:bounded:geometry}
that there is a hyperbolic manifold $Q \in \bdry B_X$ in the boundary
of the Bers slice $B_X$ so that the positive end $E^+_Q$ of $Q$ is
degenerate, and there is an orientation and marking preserving
bi-Lipschitz diffeomorphism 
$$h
\colon E^+ \to E^+_Q.$$
By Minsky's bounded geometry theorem (Theorem~\ref{Minsky:bounded}),
the bi-Lipschitz diffeomorphism type of the end $E^+_Q$ depends only
on $\nu(E^+_Q)$.  Since $E^+$ and $E^+_Q$ are bi-Lipschitz
diffeomorphic, we have $\nu(E^+) = \nu(E^+_Q)$.  
If $[\mu]$ is any limit of $\{Y_n\}$ in $\pl(S)$, it follows from
\cite[Thm. 6.1]{Brock:boundaries}
that either $\nu(E^+) = |\mu|$ or $\nu(E^-) = |\mu|$.  
The fact that $\nu(E^+) = |\mu|$ is due to Thurston (see
\cite[Sec. 9.2]{Thurston:book:GTTM}), but the argument makes delicate
use of interpolations of pleated surfaces through $Q(X_n,Y_n)$.  We
present an alternative argument in which a geodesic in the complex
of curves plays a similar role to that of the carefully chosen path in
$\ml(S)$ from Thurston's original argument.

By strong convergence, there are compact cores $\calM_n \subset
Q(X_n,Y_n)$ converging geometrically to a compact core $\calM \subset
N$ so that each compact core separates the ambient manifold into positive
and negative ends.  It follows that if $\nu(E^-) = |\mu|$ then
there are simple closed curves $\gamma_n$ with $[\gamma_n ] \to [\mu]$
for which the geodesic representative $\gamma_n^*$ of $\gamma_n$ 
in $Q(X_n,Y_n)$ always lies in the negative end of
$Q(X_n,Y_n)$ for $n$ sufficiently large.  

Let $\beta_n \in \ml(S)$ denote the bending lamination for the convex core
boundary component of $Q(X_n,Y_n)$ that faces $Y_n$ (see
\cite{Thurston:book:GTTM,Epstein:Marden:convex}).  By a theorem of
M. Bridgeman \cite[Prop. 2]{Bridgeman:bending}, there is a $K>0$ so
that the lengths $\ell_{Y_n}(\beta_n)$  are uniformly bounded, and
since $Y_n \to [\mu]$ in $\pl(S)$ there are $\mu_n \in
\mu$ with $[\mu_n ] \to [\mu]$ so that $\ell_{Y_n}(\mu_n)$ are also
uniformly bounded (see \cite[Thm. 2.2]{Thurston:hype2},
\cite[Sec. 9]{Thurston:stretch}).  It follows 
that the projective classes $[\beta_n]$ converge up to subsequence to
a limit $[\mu']$ with $i(\mu',\mu) = 0$, so we have $|\mu| = |\mu'|$.

Given $\delta > 0$, we may construct a nearly straight train track
$\tau_n$ carrying $\beta_n$ (see
\cite[Sec. 8]{Thurston:book:GTTM} or
\cite[Lem. 5.2]{Brock:length}),
so that leaves of $\beta_n$ lie within $\delta$ of $\tau_n$, as do
the geodesic representatives of sufficiently close approximates to
$\beta_n$ by simple closed curves.  Diagonalizing, 
there are simple closed curves $\eta_n$ with $[\eta_n] \to [\mu']$
so that the geodesic representatives $\eta_n^*$ of $\eta_n$ in
$Q(X_n,Y_n)$ lie in the positive ends of $Q(X_n,Y_n)$ for all $n$
sufficiently large.

Joining $\gamma_n$ to $\eta_n$ by a geodesic $g_n$ in $\calC(S)$, it
follows from \cite{Klarreich:boundary} that the entire geodesic $g_n$
converges to $|\mu|$ in $\bdry \calC(S)$. Thus, for any sequence
$\alpha_n$ of curves so that $\alpha_n$ corresponds to a vertex of
$g_n$, we have $[\alpha_n] \to [\mu]$. Since successive pairs of
curves of along $g_n$ have zero intersection, we may realize each such
pair by a pleated surface, one of which, say $X_n$, must intersect
$\calM_n$.  The sequence of pleated surfaces $X_n$ in $Q(X_n,Y_n)$ has
a limit $X_\infty$ in $N$ (after passing to a subsequence, see
\cite[Prop. 5.9]{Thurston:hype1} \cite{Canary:Epstein:Green}) realizing
$\mu$, a contradiction.

It follows that $\nu(E^+) = |\mu|$, and we conclude that the
bi-Lipschitz diffeomorphism type of the end $E^+$ depends only on the
sequence $\{Y_n\}$.

\smallskip

If $N$ has a closed geodesic $\gamma^*$ for which $\ell_N(\gamma^*) <
\ell_0/2$, then there is an $I$ so that for all $n > I$, we have
$$\ell_{Q(X_n,Y_n)} (\gamma) < \ell_0.$$ Otal's theorem
(Theorem~\ref{theorem:otal}) implies that $\gamma^*$ is the geodesic
representative of a simple closed curve $\gamma$ on $S$, and that
$\gamma^*$ is unknotted in
each $Q(X_n,Y_n)$.  It follows that the complement $Q(X_n,Y_n) \setminus
\tube^\epsilon(\gamma)$ of the Margulis tube for $\gamma$ in
$Q(X_n,Y_n)$ has the homeomorphism type
$$Q(X_n,Y_n) \setminus \tube^{\varepsilon}(\gamma) \cong S \times
\reals \setminus \calN(\gamma \times \{0\})$$
where $\calN(\gamma \times \{0\})$ denotes an embedded solid torus
neighborhood of $\gamma \times \{0\}$ in $Q(X_n,Y_n)$.

\newcommand{\collar}{{\bf collar}}

We assume for simplicity that the curve $\gamma$ is non-separating on
$S$; the separating case presents no new subtleties.  For reference
let
$$T = S \setminus \collar(\gamma)$$
be the essential subsurface in the complement of a standard open
annular collar of $\gamma$.

Applying Theorem~\ref{theorem:bilip}, we may send the cone-angle at
$\gamma$ to zero keeping the conformal boundary fixed, to obtain a
manifold $M_\gamma(X_n,Y_n)$ with a rank-2 cusp $\cusp(\gamma)$ at
$\gamma$ and a uniformly bi-Lipschitz diffeomorphism of pairs $$h_n
\colon \left( Q(X_n,Y_n) \setminus \tube^{\varepsilon}(\gamma),
\bdry \tube^{\varepsilon}(\gamma) \right ) \to 
\left( M_\gamma(X_n,Y_n) \setminus \cusp(\gamma), \bdry \cusp(\gamma)
\right)$$ so that $h$ is marking 
preserving on the ends of $Q(X_n,Y_n)$, in a sense to be made precise
presently.   

For simplicity of notation, let
$Q_n = Q(X_n,Y_n)$
and
$M_n = M_\gamma(X_n,Y_n).$  
We pause to elaborate briefly on the structure of $M_n$.
Since the geodesic $\gamma^*$ is unknotted, we may choose a product
structure $F_n \colon S \times \reals \to Q_n$ so that
$\gamma^* = F_n( \gamma \times \{0 \})$ and so that there is a
standard
tubular neighborhood 
$$V_\gamma = \collar(\gamma) \times (-1/2,1/2)$$ of 
$\gamma \times \{0\}$ in $S\times \reals$ with $F_n(V_\gamma) =
\tube^{\varepsilon}(\gamma)$.  

We take $\calM_n = F(S \times [-1,1])$ as a compact core for $Q_n$,
and we denote by $E^+_n = F(S \times (1,\infty))$ and $E^-_n = F(S
\times (-\infty, -1))$ the positive and negative ends of $Q_n$ with
respect to $\calM_n$.  The diffeomorphism $h_n$ sends $\calM_n
\setminus \tube^\varepsilon(\gamma)$ to a relative compact core for $M_n$ and
sends $E^+_n$ and $E^-_n$ to ends of $M_n$.  We will refer to
$h_n(E^+_n)$ and $h_n(E^-_n)$ as the {\em positive} and {\em negative}
ends of $M_n$.  The positive and negative ends of $M_n$ are implicitly
marked by the compositions $$h_n \compos F(S \times \{1\}) \ \ \
\text{and} \ \ \ h_n \compos F(S \times \{-1\}).$$ 
It follows from the fact that the cone-deformation preserves the
conformal boundary that in 
this marking
the positive ends of $Q_n$ and $M_n$ are compactified by the same surface
$Y_n \in \Teich(S)$ and the negative ends of $Q_n$ and $M_n$ are
compactified by the same surface $X_n \in \Teich(S)$.

Let $\omega_n \in Q_n \setminus \tube^{\varepsilon}(\gamma)$ be a
choice of base-frame for $Q_n$ so that the lifts $(Q_n,\omega_n)$ to
$AH_\omega(S)$ converge strongly to $(N,\omega)$.  Let $\omega_n' =
h_n(\omega_n)$ be base-frames for $M_n$.  By Gromov's compactness
theorem, we may pass to a subsequence and extract a geometric limit
$(h_\infty,(Q_\infty,\omega_\infty),(M_\infty,\omega_\infty'))$ of the
triples $(h_n,(Q_n,\omega_n),(M_n,\omega_n'))$ where $h_\infty$ is a
uniformly bi-Lipschitz diffeomorphism between smooth submanifolds of
$Q_\infty$ and $M_\infty$ obtained as a limit of the uniformly
bi-Lipschitz diffeomorphisms $h_n$.

Since the convergence of $Q_n \to N$ is strong, we have $Q_\infty =
N$, and $h_\infty$ has domain $N \setminus
\tube^{\varepsilon}(\gamma)$. 
By the above discussion, the
mappings $h_n$ induce the markings on the positive ends of $M_n$, so
$h_\infty$ gives a marking and orientation preserving uniformly
bi-Lipschitz diffeomorphism of the positive end $E^+$ of $N$ with the
positive end $E^+_M$ of $M_\infty$. 

The covering $Q_\gamma$ of $M_\infty$ corresponding to $\pi_1(E^+_M)$
lies in $AH(S)$ and has the following properties:
\begin{enumerate}
\item $Q_\gamma$  has positive end isometric to $E^+_M$, 
\item $Q_\gamma$ has cuspidal thin part $P_\gamma \subset Q_\gamma$, a rank
one cusp whose inclusion on $\pi_1$ is conjugate to $\langle \gamma \rangle$, 
\item the pared submanifold $Q_\gamma \setminus P_\gamma$ has a negative end
$$E_\gamma^- \cong T \times (-\infty,
0].$$
\end{enumerate}
We claim that the end $E_\gamma^-$ of the pared submanifold $Q_\gamma
\setminus P_\gamma$ is geometrically finite. 

Since it suffices to show that the cover $\widetilde{Q_\gamma}$ of
$Q_\gamma$ corresponding to $\pi_1(E_\gamma^-) \cong \pi_1(T)$ is
quasi-Fuchsian, assume otherwise.  Then there is some degenerate end
$E$ of the {\em pared submanifold} of $\wt Q_\gamma$, namely, the
complement $\wt Q_\gamma \setminus \wt P_\gamma$ of the cuspidal thin
part of $\wt Q_\gamma$. 
Since the manifold $Q_\gamma$ is not a
surface bundle over the circle, the covering theorem (see
\cite{Thurston:book:GTTM}, \cite{Canary:inj:radius}) implies that the
end $E$ has a neighborhood $U$ that covers a degenerate end of
$M_\infty$ finite-to-one.  But the manifold $M_\infty$ has exactly two
degenerate ends, each homeomorphic to a product with the larger
surface $S$ with a half-line.

It follows that the negative end of $Q_\gamma$ is geometrically
finite.  Let $Z \in \Teich(T)$ denote the associated conformal
boundary component.   Choose a surface $X \in \Teich(S)$ so that
$\ell_X(\gamma) < \ell_0/4$.  By a theorem of Bers (see 
\cite[Thm. 3]{Bers:bdry}  or \cite[Prop. 6.4]{McMullen:iter}),
then, the manifolds
$Q(X,Y_n)$ have the property that
$$\ell_{Q(X,Y_n)}(\gamma) < \ell_0/2$$ for all $n$.

Performing the same process as above with $X_n = X$, we arrive
at a manifold $Q_\gamma'$ whose positive end $\check E^+$ is
bi-Lipschitz diffeomorphic to the positive end of 
the limit $Q_\infty$ of $Q(X,Y_n)$ after passing to a subsequence.
The negative end of $Q_\gamma'$ has conformal boundary surface $Z'
\in \Teich(T)$.  

The manifold $Q_\gamma$ is the cover associated to the positive end of
$$\lim_{n \to \infty} M_\gamma(X_n,Y_n) = M_\infty$$  and the manifold
$Q_\gamma'$ 
is the cover associated to the positive end of the limit $$\lim_{n \to
\infty} M_\gamma(X,Y_n).$$
But letting $Q_\gamma(n)$ be the cover of $M_\gamma(X_n,Y_n)$
associated
to $Y_n$ and letting $Q_\gamma(n)'$ be the cover of $M_\gamma(X,Y_n)$
associated to $Y_n$, we have
$$\bdry Q_\gamma(n) = Z_n \disjunion Y_n \ \ \ \text{and} 
\ \ \ \bdry Q_\gamma(n)' = Z_n' \disjunion Y_n$$
where $Z_n$ converges to $Z$ in $\Teich(T)$ and $Z_n'$ converges to
$Z'$ in $\Teich(T)$.  Letting $K_n = d_T(Z_n, Z_n')$ be the
Teichm\"uller distance from $Z_n$ to $Z_n'$, the manifolds
$Q_\gamma(n)$ and $Q_\gamma(n)'$ have quasi-isometric distance
$$d_{\rm qi}(Q_\gamma(n),Q_\gamma(n)') < K_n'$$ where $K_n'$ depends
only on $K_n$ (see e.g. \cite[Thm 2.5]{McMullen:book:RTM}).  By lower
semi-continuity of the quasi-isometric distance on $AH(S)$
\cite[Prop. 3.1]{McMullen:book:RTM}
we have a marking preserving bi-Lipschitz diffeomorphism
$$ \Phi \colon Q_\gamma \to Q_\gamma'.$$

Restricting the diffeomorphism $\Phi$ to the positive ends of
$Q_\gamma$ and $Q_\gamma'$, it follows that the positive end of $N$ is
bi-Lipschitz diffeomorphic to the positive end of $$Q_\infty = \lim_{n
\to \infty} Q(X,Y_n).$$ Since for any surface $X'\in \Teich(S)$, the
manifolds in the sequences $Q(X,Y_n)$ and $Q(X',Y_n)$ are
uniformly bi-Lipschitz diffeomorphic, the bi-Lipschitz diffeomorphism
type of $E^+$ does not depend on $X$ (by another application of
\cite[Prop. 3.1]{McMullen:book:RTM}), so the theorem follows.
\qed

\section{Proof of the main theorem}
\label{section:proof}

In this section, we assemble our results to give the proof of our main
approximation theorem (Theorem~\ref{theorem:main}).  The proof
naturally breaks into cases based on the homotopy type of the
manifold we wish to approximate.  We first treat the following case.
\begin{theorem}
Let $M \in AH(S)$ be cusp-free.  Then $M$ is an algebraic limit of
quasi-Fuchsian manifolds.
\label{theorem:surface} 
\end{theorem}

\bold{Proof:} 
By Theorem~\ref{theorem:ends:realized} (or  the main theorem of \cite{Bromberg:bers}), it suffices to consider
the case when $M$ is doubly degenerate.  Let $F \colon S \times \reals
\to M$ be a smooth product structure on $M$ and let $\calM = F(S \times
 [-1,1])$ denote a compact core for $M$.  Let 
$$E^-  = F(S \times (-\infty, -1))  \ \ \ \text{and} \ \ \ 
E^+ = F(S \times (1,\infty))$$ 
denote the positive and negative ends of $M$.

By Theorem~\ref{theorem:ends:realized}, the ends $E^+$ and $E^-$ are
realizable by manifolds 
$$Q^+ = \lim_{n \to \infty} Q(X,Y_n) \ \ \
\text{and} \ \ \ Q^- \lim_{n \to \infty} Q(X_n,Y)$$ in the Bers
boundaries $\bdry B_X^+$ and $\bdry B_Y^-$.

Moreover, if $\nu^+$ and $\nu^-$ are the end-invariants for $M$,
Corollary~\ref{corollary:pml} guarantees that after passing to a
subsequence we may assume that we have the convergence
$$X_n \to [\mu^-] \ \ \ \text{and} \ \ \ Y_n \to [\mu^+]$$
in Thurston's compactification $\pl(S)$, and that the support of
$\mu^-$ and $\mu^+$ is given by
$$|\mu^-| = \nu^- \ \ \ \text{and} \ \ \ |\mu^+| = \nu^+.$$

Since $\nu^-$ and $\nu^+$ are the end-invariants of the doubly
degenerate manifold $M$, it follows that $\mu^-$ and $\mu^+$ {\em
bind the surface}: for any simple closed curve $\alpha$, we have
$$i(\alpha,\mu^-) + i(\alpha,\mu^+) >0.$$

Thus, by Thurston's {\em double limit theorem}
(\cite[Thm. 4.1]{Thurston:hype2}, see also \cite{Otal:book:fibered}), 
the sequence
$$\{Q_n\}_{n=1}^\infty = \{Q(X_n,Y_n)\}_{n=1}^\infty$$
converges after passing to a subsequence.

By an application of the continuity
of the length function \cite[Thm. 2]{Brock:length} we may pass to a
subsequence so that the sequence $Q(X_n,Y_n)$
converges to a limit $M' \in AH(S)$ for which $\nu^+$ and $\nu^-$ are
the end-invariants of $M'$.  It follows that the limit $M'$ is doubly 
degenerate.

Applying Theorem~\ref{theorem:isolation}, the ends of $M'$ are
geometrically isolated: the positive end of $M'$ is bi-Lipschitz
diffeomorphic to the positive end of $Q^+$ and the negative end of
$M'$ is bi-Lipschitz diffeomorphic to the negative end of $Q^-$.
Since $Q^+$ and $Q^-$ realize the ends of $M$, it follows 
that there is a smooth product structure 
$$F' \colon S \times \reals \to M'$$ so that
\begin{enumerate}
\item $F'$ decomposes $M'$ into a compact core $\calM' = F(S \times [-1,1])$ and ends 
$(E^-)' = F(S \times (-\infty,-1))$ and
$(E^+)' = F(S \times (1,\infty))$, 
\item there are marking preserving bi-Lipschitz diffeomorphisms
$$h^+ \colon E^+ \to (E^+)' \ \ \ \text{and}
\ \ \ h^- \colon E^- \to (E^-)'$$ 
between positive and negative ends
of $M$ and $M'$ so that 
$$h^+ (F(x,t)) = F'(x,t) \ \ \ \text{ and } \ \ \ 
h^- (F(x,t)) = F'(x,t).$$  
\end{enumerate}
By compactness of $\calM$, the extension 
$$ h \colon M \to M'$$ of 
$h^+$ and $h^-$ across $\calM$
defined by setting $h(F(x,t)) = F'(x,t)$
is a single marking preserving bi-Lipschitz diffeomorphism from $M$ to
$M'$.
Applying Sullivan's rigidity theorem \cite{Sullivan:linefield}, it
follows that $h$ is homotopic to an isometry, and we may conclude
that $Q(X_n,Y_n)$ converges to $M$.
\qed

To complete the proof of Theorem~\ref{theorem:main}, we now treat the
case when $M$ is a general infinite volume complete hyperbolic
3-manifold with no cusps incompressible ends.  The theorem in this
case essentially follows directly from the case when $M$ has the
homotopy type of a surface; there are two issues to which we alert the
reader:
\begin{enumerate}
\item In the surface case, we applied Thurston's  double limit theorem 
to show the 
approximates converge up to subsequence.
Here, an analogous compactness result is necessary
(\cite[Thm. 2.4]{Ohshika:compact},
cf. \cite{Thurston:hype3}).
\item  Examples  constructed by  Anderson and Canary (see
\cite{Anderson:Canary:pages}) illustrate that homeomorphism-type
need not persist under algebraic limits of hyperbolic 3-manifolds.
Since we control the ends of $M_n$, and thence the  peripheral
structure of $\pi_1(M_n)$, we may prevent such a topological cataclysm by an
application Waldhausen's theorem
\cite{Waldhausen:irreducible}
\cite[Thm. 13.7]{Hempel:book}.
\end{enumerate}
\begin{theorem}  
Let $M$ be an infinite volume complete hyperbolic 3-manifold with
incompressible ends and no cusps.  Then $M$ is an algebraic limit 
of geometrically finite manifolds.
\label{theorem:mixed}
\end{theorem}

\bold{Proof:}  
By Theorem~\ref{theorem:surface}, we may assume that
$M$ is not homotopy equivalent to a surface.
By Bonahon's theorem, there is a compact 3-manifold $N$ so that 
$$M \cong \interior(N).$$  Thus, $M$ lies in $AH_0(N)$, the subset of
$AH(N)$ consisting of marked hyperbolic 3-manifolds $(f \colon N \to
M)$ for which $f \vert_{\interior(N)}$ is homotopic to a homeomorphism.

As in the outline, we approximate the tame manifold $M$
end-by-end, and combine the approximations into one sequence of
geometrically finite hyperbolic manifolds $M_n$ homeomorphic to $M$
that converge to $M$.

Let $\calM$ be a compact core for $M$ and 
let $E$ be a degenerate end of $M \setminus \calM$.
Let $Q$ be the cover of $M$ corresponding to $\pi_1(E)$.
Then $Q$ lies in $AH(S)$, 
where $S = \closure{E} \cap \calM$.  
Assume
$M$ is oriented so that $E$ lifts to the positive end $\wt E$ of $Q$ in the
cover.

By the covering theorem (see \cite[Ch. 9]{Thurston:book:GTTM} or
\cite[Main Thm.]{Canary:inj:radius}), 
we claim the manifold $Q$ is a singly degenerate manifold with no
cusps and degenerate end $\wt E$.  To see this, note first that $M$
has no cusps, so neither does $Q$.  The only alternative is then that
$Q$ is doubly degenerate which implies that the covering $Q\to M$ is
finite-to-one, since $M$ is not a surface bundle over the circle.  
But
if $Q$ covers $M$ finite-to-one, the manifold
$M$ is itself homotopy equivalent to a surface finitely covered by $S$.

The geometrically finite locus $GF_0(N) = AH_0(N)$ consists of
geometrically finite hyperbolic 3-manifolds $M$ homeomorphic to
$\interior(N)$.  Realizing the degenerate ends of $M$ in the
appropriate Bers boundaries, we will obtain surfaces that determine
candidate approximates for $M$ in the interior of $AH_0(N)$.  The
interior of $AH_0(N)$ is typically denoted by $MP_0(N)$, the {\em
minimally parabolic} structures on $\interior(N)$, namely, structures
on $\interior(N)$ with only finite volume cusps.  Since $N$ is assumed
to have no torus boundary components, $MP_0(N)$ is simply the
cusp-free hyperbolic structures on $\interior(N)$ in the case at hand.

Choosing a reference compact core $\calM$ for $M$,
let $\calE_1, \ldots, \calE_{p}$ denote the geometrically finite ends
of $M$ and let $E_1, \ldots, E_{q}$ 
denote its simply degenerate ends.
Let $S_j \subset \bdry N$ denote the boundary component of $N$ lying
in the closure of the end $\calE_j$, $j = 1,\ldots, {p}$ and let $T_k
\subset \bdry N$ denote the boundary component of $N$ lying in the
closure of the end $E_k$ where $k = 1, \ldots, q$.  Then the space
$MP_0(N)$ admits the parameterization $$MP_0(N) =\Teich(\bdry N) =
\prod_{j=1}^{p} \Teich(S_j) \times
\prod_{k=1}^{q}\Teich(T_k).$$  (see
\cite[Thm. 14]{Kra:def}, \cite[Thm. 1.3]{Thurston:hype2}, or
\cite{McMullen:iter}).   

Let $\{Y_k(n)\}_{n=1}^\infty \subset \Teich(T_k)$ denote the sequences
of surfaces obtained from Theorem~\ref{theorem:ends:realized} so
that the limit $$Q_k =
\lim_{n \to \infty} Q(Y_k(0),Y_k(n))$$ realizes the end $E_k$, $k = 1,
\ldots, p$.  Let $$h_k \colon \hat E_k \to E_k$$ be the marking
preserving bi-Lipschitz diffeomorphism from the positive end of $Q_k$
to $E_k$ coming from Theorem~\ref{theorem:ends:realized}. Let $X_j$
be the conformal boundary component compactifying $\calE_j$, for $j =
1,\ldots, {p}$ and let $$g_j \colon \hat \calE_j \to \calE_j$$ be the
marking preserving bi-Lipschitz diffeomorphism from the end of the
Fuchsian manifold $Q(X_j, X_j)$.

In the above parametrization, we let $M_n$ be determined by $$M_n =
(X_1,\ldots, X_{p},Y_1(n),\ldots,Y_p(n)) \in \Teich(\bdry N).$$ Let $f_n
\colon N \to M_n$ be the implicit homotopy equivalences marking the
manifolds $M_n$.  We claim that the sequence $M_n$ converges up to
subsequence in $AH(N)$.  

Corollary~\ref{corollary:pml} guarantees $Y_k(n)$ may be chosen to
converge to a limit $[\mu_k] \in \pl(S)$ with support of $|\mu_k| =
\nu_k$ where $\nu_k = \nu(E_k)$.  Thus, there are simple closed curves
$\gamma_n$ on $T_k$ and positive real weights $t_n$ so that $t_n
\gamma_n$ converges in $\ml(T_k)$ to $\mu_k$, and the lengths
$\ell_{Y_k(n)}(t_n \gamma_n)$ remain bounded. Applying a Theorem of
K. Ohshika generalizing Thurston's compactness theorem (see
\cite[Thm. 2.4]{Ohshika:compact},
\cite{Thurston:hype3})  we conclude that
the sequence $M_n$ converges after passing to a subsequence.
We pass to a subsequence so that $M_n$ converges algebraically to
$M_\infty$, and geometrically to a manifold $M_G$ 
covered by $M_\infty$ by a local
isometry  $\pi \colon M_\infty \to M_G$. 
Let $f_\infty \colon N \to
M_\infty$ denote the marking on $M_\infty$. 

We will exhibit a compact core $\calM_\infty$ for $M_\infty$ so that
any homotopy equivalence $f_\infty' \colon N \to \calM_\infty$
homotopic to $f_\infty$ is homotopic to a homeomorphism, and the ends
$M_\infty \setminus
\calM_\infty$ are exactly
the ends $\hat E_k$ and $\hat \calE_j$, where $k = 1, \ldots, p$ and
$j = 1,\ldots {p}$.

First consider the covers $\calQ_j(n)$ of $M_n$ corresponding to
$\pi_1(X_j)$.  The manifolds $\calQ_j(n)$ range in the Bers slice $B_{X_j}$,
and the surface $X_j$ persists as a conformal boundary component of
the algebraic and geometric limits $\calQ_j(\infty)$ and $\calQ_j^G$ of
$\{\calQ_j(n)\}$ 
(after potentially passing to a further subsequence; see,
e.g. \cite[Prop. 2.3]{Kerckhoff:Thurston}).  

The end $\calE_j(n)$ of $\calQ_j(n)$ cut off by the boundary of a smooth
neighborhood of the convex core boundary component facing $X_j$ embeds
in the covering projection to $M_n$.  
The ends $\calE_j(n)$ converge
geometrically to the geometrically finite end $\calE_j(\infty)$ of the
geometric limit $\calQ_j^G$, and this end is compactified by $X_j$.  
It follows that there are marking preserving smooth embeddings
$\varphi_j^n$ of $\calE_j(\infty)$ into $M_n$ that converge $C^\infty$ to an
isometric embedding $\varphi_j \colon \calE_j(\infty) \to M_G$.  

The isometric embedding $\varphi_j$ is 
marking preserving in that 
$$(\varphi_j)_* \vert_{(\pi_1(S_j))} = \pi_* \compos
(f_\infty)_* \compos (\iota_j)_*$$ 
where $\iota_j \colon S_j \to N$ is the inclusion map.
Thus, the end $\calE_j(\infty)$ lifts to an end 
$\bar \calE_j$ of the algebraic limit
$M_\infty$.

Consider, on the other hand, covers $Q_k(n)$ of $M_n$ corresponding to
$Y_k(n)$.  The manifolds $Q_k(n) \subset AH(S)$ converge algebraically to
$Q_k(\infty)$ and there is a measured lamination $\mu_k$ so that
$\ell_{Q_k(n)}(\mu_k) \to 0$.  The measured lamination $\mu_k$ ``fills'' the
surface $T_k$: for any essential simple closed curve
$\gamma$ on $T_k$, we have $i(\gamma,\mu_k) \not=0$.  Applying
\cite[Thm. 2]{Brock:length}, the support $|\mu_k| = \nu_k$ is an ending
lamination for $Q_k(\infty)$, which implies that $Q_k(\infty)$ has a
degenerate end $E_k(\infty)$ with ending lamination $\nu_k$.  By an
argument using the covering theorem and
\cite[Lem. 3.6]{Jorgensen:Marden:geomconv}, the end $E_k(\infty)$ embeds
in the geometric limit $Q_k^G$ (see 
\cite[Prop. 5.2]{Anderson:Canary:cores}).  Just as above, then, we have a limiting 
isometric embedding $\phi_k \colon E_k(\infty)
\to M_\infty$ 
which is marking preserving in the sense that $$(\phi_k)_*
\vert_{(\pi_1(T_k))} = \pi_* \compos (f_\infty)_* \compos (i_k)_*$$
where $i_k \colon T_k \to N$ is the inclusion map.  Let $\bar E_k$
denote this end of $M_\infty$.

By an application of Waldhausen's theorem (see
\cite{Waldhausen:irreducible} \cite[Thm. 13.7]{Hempel:book}), 
the homotopy equivalence $f_\infty \colon N \to M_\infty$ is homotopic to a
homeomorphism to a compact core $\calM_\infty$ for $M_\infty$ that
cuts off the geometrically finite ends $\bar \calE_j$ and the simply
degenerate ends $\bar E_k$.

We claim the ends $\bar \calE_j$ and $\bar E_k$ have no cusps.
Since $\bar \calE_j$ is a geometrically finite end compactified by the
closed surface $X_j$ there is no isotopy class in $\bar \calE_j$ with
arbitrarily short length.  Likewise, $\bar E_k$ is a degenerate end
with ending lamination $\nu_k = |\mu_k|$ that fills the surface $T_k$.   By
results of Thurston and Bonahon
(see \cite[Ch. 9]{Thurston:book:GTTM}, \cite[Prop. 3.4]{Bonahon:tame})
if $\gamma \subset T_k$ is a simple closed curve represented by a cusp
in $\bar E_k$ ($\gamma$ has representatives 
with arbitrarily
short representatives in the end $\bar E_k$) then we have 
$$i(\gamma,\mu_k) = 0,$$  contradicting the fact that
$\nu_k$
fills the surface.

It follows that $M_\infty$ is cusp-free, since any element $g \in
\pi_1(N)$ for which $f_\infty(g)$ is parabolic must have arbitrarily
short representatives in its free homotopy class exiting some end
$\bar \calE_j$ or $\bar E_k$ of $M_\infty$.  By a theorem of Anderson
and Canary
\cite[Cor. G]{Anderson:Canary:cores} we may conclude 
that $M_n$ converges strongly to $M_\infty$.

It follows that each sequence of covers $\calQ_j(n)$ or $Q_k(n)$ converges
strongly to a 
limit in a Bers boundary; otherwise the limit would be doubly
degenerate, which is ruled out once again by the covering theorem.  

  It follows that the geometric limit of the quasi-Fuchsian
manifolds $\calQ_j(n)$ or $Q_k(n)$ 
is uniformly bi-Lipschitz diffeomorphic to the
quasi-Fuchsian manifolds (either $Q(X_j,X_j)$ or $Q(Y_k(0),Y_k(n))$)
whose limit realizes
the corresponding end
$\calE_j$ or  $E_k$.  Thus each end of $M_\infty$
admits a marking preserving bi-Lipschitz diffeomorphism to the end of
the corresponding realization in a Bers boundary.

As these ends admit marking-preserving bi-Lipschitz diffeomorphisms to
the ends of $M$, we may extend the corresponding bi-Lipschitz
diffeomorphisms across the compact cores to obtain a single marking
preserving bi-Lipschitz diffeomorphism $$\Psi \colon M_\infty \to M.$$
Since the corresponding conformal structures compactifying the
geometrically finite ends of $M_\infty$ and $M$ are the same, we may
apply Sullivan's theorem \cite{Sullivan:linefield} to conclude that
$\Psi$ is homotopic to an isometry, and we have $$M = M_\infty =
\lim_{n \to \infty} M_n.$$ As the manifolds $M_n$ are geometrically
finite, the proof is complete.
\qed

\bibliographystyle{math}
\bibliography{math}

\noindent{\sc \scriptsize Department of Mathematics, University of Chicago,
5734 S. University Ave.,
Chicago, IL  60637}

\noindent{\sc \scriptsize Department of Mathematics,
California Institute of Technology,
253-37,
Pasadena, CA 91125} 

\end{document}